%% file: main.tex
\documentclass[journal]{IEEEtran}

\usepackage{cite}

\usepackage[cmex10]{amsmath}

\usepackage[T1]{fontenc}
\usepackage{txfonts}

\usepackage{mathtools} 

\usepackage{accents}
\makeatletter
\def\widebar{\accentset{{\cc@style\underline{\mskip11mu}}}}
\makeatother

\usepackage{array}

\ifCLASSOPTIONcompsoc
  \usepackage[caption=false,font=normalsize,labelfont=sf,textfont=sf]{subfig}
\else
  \usepackage[caption=false,font=footnotesize]{subfig}
\fi

\usepackage{booktabs}
\usepackage{tabularx}

\usepackage{url}
\usepackage{multirow}

\usepackage{subfiles}

\newlength{\nomenlabelindent}
\newenvironment{nomenclature}{%
\newcommand\entry[2]{%
   \hangindent\nomenlabelindent\noindent\makebox[\nomenlabelindent][l]{##1\quad}\ignorespaces##2\par}%
   }

\begin{document}

	\title{Hydrogen Supply Chain Planning with Flexible Transmission and Storage Scheduling}
	%

% \begingroup
% \author{Guannan He\footnotemark[1]~~, 
% 		Dharik S. Mallapragada\footnotemark[1]
% % 		~~\footnotemark[3]
% 		~~,
% 		Abhishek Bose\footnotemark[1]~~,
% 		Clara F. Heuberger\footnotemark[2]~~,  
% 		Emre Gençer\footnotemark\footnotemark[1]
% % 		~~\footnotemark[3]
% 	}
% \endgroup

	\author{Guannan~He,
% 	~\IEEEmembership{Member,~IEEE,}
		Dharik S. Mallapragada,
		Abhishek Bose,
		Clara F. Heuberger,
		and~Emre Gen\c{c}er% <-this % stops a space
		\thanks{This work was partially supported by Shell New Energies Research and Technology, Amsterdam, Netherlands, and the Low-Carbon Energy Centers at MIT Energy Initiative. Guannan He, Dharik S. Mallapragada, Abhishek Bose, and Emre Gençer are with MIT Energy Initiative, Massachusetts Institute of Technology, Cambridge, MA 02139, USA. Clara F. Heuberger is with Shell Global Solutions International B.V., Shell Technology Centre Amsterdam, 1031 HW Amsterdam, Netherlands. (e-mail: gnhe@mit.edu)}% <-this % stops a space

		\vspace{-2em}}

	%\markboth{IEEE Transactions on Smart Grid}%
	%{}
	
	\maketitle
	
    \begin{abstract}
    Hydrogen is becoming an increasingly appealing energy carrier, as the costs of renewable energy generation and water electrolysis continue to decline. Developing modelling and decision tools for the H$_{2}$ supply chain that fully capture the flexibility of various resources is essential to understanding the overall cost-competitiveness of H$_{2}$ use. To address this need, we have developed a H$_{2}$ supply chain planning model that determines the least-cost mix of H$_{2}$ generation, storage, transmission, and compression facilities to meet H$_{2}$ demands and is coupled with power systems through electricity prices. We incorporate flexible scheduling for H$_{2}$ trucks and pipeline, allowing them to serve as both H$_{2}$ transmission and storage resources to shift H$_{2}$ demand/production across space and time. The case study results in the U.S. Northeast indicate that the proposed framework for flexible scheduling of H$_{2}$ transmission and storage resources is critical not only to cost minimization but also to the choice of H$_{2}$ production pathways between electrolyzer and centralized natural-gas-based production facilities. Trucks as mobile storage could make electrolyzer more competitive by providing extra spatiotemporal flexibility to respond to the electricity price variability while meeting H$_{2}$ demands. The proposed model also provides a reasonable trade-off between modeling accuracy and computational time.
    \end{abstract}
	
	\begin{IEEEkeywords}
        Hydrogen supply chain planning, electrolytic hydrogen production, mobile storage, flexibility
        \vspace{-1em}
	\end{IEEEkeywords}

	\IEEEpeerreviewmaketitle
	
    % \bstctlcite{IEEEexample:BSTcontrol}	
    \subfile{Sections/0_nomenclature}

    \subfile{Sections/1_introduction}

    % \subfile{Sections/literature}
    
    \subfile{Sections/2_flexibility_modelling}
    
    \subfile{Sections/3_model_formulation}
    
    % \subfile{Sections/data_setup}
    
    \subfile{Sections/4_results}

    \subfile{Sections/5_conclusion}

	\subfile{Sections/Acknowledgment}

% 	\subfile{Sections/Appendix}
% 	\subfile{Sections/Brief_proof}

	\ifCLASSOPTIONcaptionsoff
	\newpage
	\fi

	
	
	
	\bibliography{IEEEabrv,MyLibrary}
% 	\bibliography{MyLibrary}
% 	\bibliographystyle{IEEEtran}
	% the LaTeX parser from getting confused when it sbattery the complicated
	
	%\end{IEEEbiography}
	
	
	%\newpage
	
	
	
	%\vfill
	

\end{document}

\title{Hydrogen Supply Chain Planning with Flexible Transmission and Storage Scheduling: Supporting Information}
\author{Guannan~He,
		Dharik S. Mallapragada,
		Abhishek Bose,
		Clara F. Heuberger,\\
		and~Emre Gençer% <-this % stops a space
		}
\date{}
\thispagestyle{empty}
\pagestyle{empty}
\maketitle

\newcolumntype{P}[1]{>{\centering\arraybackslash}p{#1}}

%\title{Data Tables for Base Scenario}
\section{Transportation Demand}

\subsection{Average Hydrogen Refuelling Demand}

% Table generated by Excel2LaTeX from sheet 'Sheet2'
\begin{table}[htbp]
  \centering
  \caption{Average Hourly Hydrogen Demand for Light-Duty (LDV) and Heavy-Duty (HDV) vehicles \cite{faf,EPA_freight}} 
  
    \begin{tabular}{ccc}
    \toprule
          & \textbf{LDV Demand} & \textbf{HDV Demand} \\
          & tonne/hour & tonne/hour \\
    \midrule
    Zone 1 & 25    & 6 \\
    Zone 2 & 159   & 33 \\
    Zone 3 & 57    & 12 \\
    Zone 4 & 123   & 21 \\
    Zone 5 & 39    & 9 \\
    Zone 6 & 55    & 46 \\
    Zone 7 & 0     & 0 \\
    \bottomrule
    \end{tabular}%
  \label{tab:Trans_Avg_demand}%
\end{table}%

\subsection{Hydrogen Refuelling Profile}
\begin{figure}[H]
    \centering
    \includegraphics[width = \textwidth]{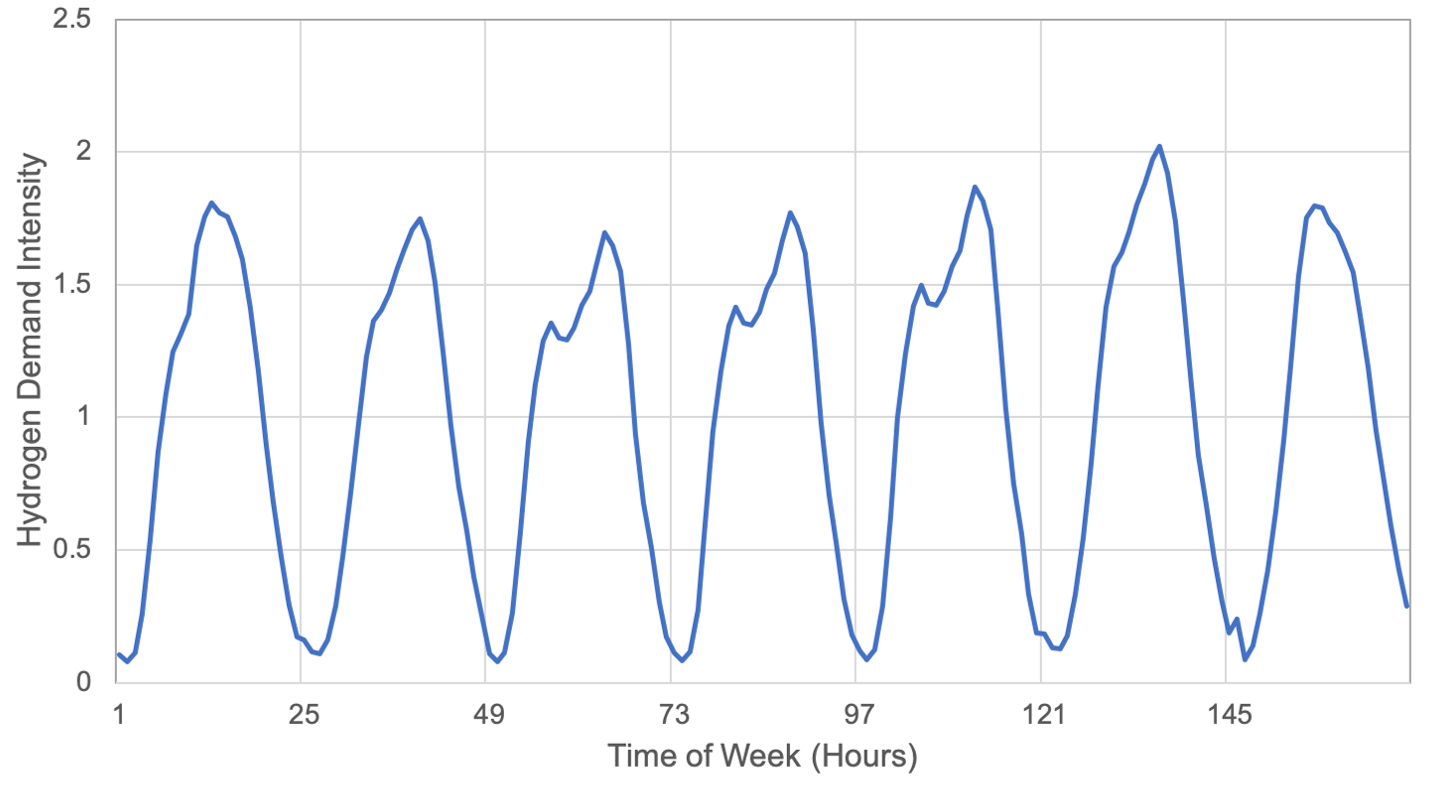}
    \caption{Hourly hydrogen refuelling profile normalized based on the mean}
    \label{fig:Unit Refueling Profile}
\end{figure}

\subsection{Hydrogen Demand Profile Calculation}
$\text{H}^{\text{D}}_{z}$ : the sum of LDV and HDV demands, as shown in Table. \ref{tab:Trans_Avg_demand}

$U_{z,t}$ : hydrogen refuelling profile, as shown in Fig. \ref{fig:Unit Refueling Profile}
\begin{equation}
\text{D}_{z,t} = \text{H}^{\text{D}}_{z} \cdot \text{U}_{z,t}
\end{equation}

% \clearpage
\section{Inter Zone Distances}

\begin{table}[H]
  \centering
  \caption{Inter-Zone Distances for Hydrogen Supply Chain Model ($L_{z->z'}$, mile)}
  \vspace{1em}
    \begin{tabular}{c|cccccc}
    \toprule
    \textbf{Zone} & \textbf{1} & \textbf{2} & \textbf{3} & \textbf{4} & \textbf{5} & \textbf{6} \\
    \midrule
    \textbf{1} & 0     & 317   & 504   & 602   & 487   & 608 \\
    \textbf{2} & 317   & 0     & 199   & 297   & 179   & 340 \\
    \textbf{3} & 504   & 199   & 0     & 99    & 158   & 333 \\
    \textbf{4} & 602   & 297   & 99    & 0     & 216   & 358 \\
    \textbf{5} & 487   & 179   & 158   & 216   & 0     & 186 \\
    \textbf{6} & 608   & 340   & 333   & 358   & 186   & 0 \\
    \bottomrule
    \end{tabular}%
  \label{tab:zone_dist}%
\end{table}%

\clearpage
\section{Electricity Price Profiles}
\textbf{Description:}
\begin{enumerate}
    \item Power system capacity expansion model used to generate prices for each zone and hour of the year (\url{http://energy.mit.edu/eps/tools/genx/})
    \item Greenfield setup considered for capacity expansion. Existing capacity considered for hydro, pumped storage and transmission infrastructure
    \item 20 representative weeks are modelled for annual operations, based on methods described in \cite{MALLAPRAGADA2020115390}
    \item Capacity expansion in transmission of electricity is not allowed in the base case model
\end{enumerate}
\textbf{Load Data}: \\
2018 NREL electrification futures study (EFS) projected load for 2050 with assumed BAU technology advancement and reference electrification.
Reference: \url{https://data.nrel.gov/files/126/EFSLoadProfile_Reference_Moderate.zip}
\\
\newline
\textbf{ Resource Characterization}:
\begin{enumerate}
    \item \textbf{Renewable Data Source}: NREL Wind Toolkit (onshore and offshore wind) and National Solar Radiation database (solar)
    \item Offshore wind is included with no capacity limits and single resource profile for zone 2 and zone 4 based on sampling sites from the NREL Wind Toolkit that overlaps with the areas for which leases have been auctioned
    \item Three types of hydropower resources included in the model: Reservoir hydro (with flexible generation and constraints on weekly energy generation limits), run-of-river hydro (inflexible generation) and Canadian hydro imports (Zone 7, flexibile with weekly energy generation limits)
\end{enumerate}

\subsection{Power System Technology Cost and Performance Assumptions}

\begin{table}[H]
  \centering
  \caption{Technology Cost Details, in 2019 US\$ }
  \resizebox{\textwidth}{!}
  {
    \begin{tabular}{lllllll}
    \toprule
    \multicolumn{1}{p{9.5em}}{\textbf{Technology}} & \multicolumn{1}{p{6em}}{\textbf{Capital Cost}} & \multicolumn{1}{p{6em}}{\textbf{Capital Cost}} & \multicolumn{1}{p{6em}}{\textbf{Lifetime}} & \multicolumn{1}{p{6em}}{\textbf{Fixed O\&M Cost}} & \multicolumn{1}{p{6em}}{\textbf{Fixed O\&M Cost}} & \multicolumn{1}{p{6em}}{\textbf{Variable O\&M Cost}} \\
          & \textit{\$/MW} & \textit{\$/MWh} & \textit{year} & \multicolumn{1}{c}{\textit{\$/(MW-year)}} & \textit{\$/(MWh-year)} & \textit{(\$/MWh)} \\
    \midrule
    Onshore Wind &          1,085,886  & -     & 30    &                       34,568  & -     & 0 \\
    Offshore Wind &          1,901,981  & -     & 30    & 48,215 & -     & 0 \\
    Utility-Scale Solar &             724,723  & -     & 30    &                       11,153  & -     & 0 \\
    Li-ion Battery &             120,540  & 125870 & 15    &                         2,452  & 3071  & 0.03 \\
    Pumped Hydro &          1,966,151  & -     & 50    &                       41,000  & -     & 0.02 \\
    CCGT  &             816,971  & -     & 30    &                       10,560  & -     & 2.77 \\
    OCGT  &             815,863  & -     & 30    &                       12,230  & -     & 7.14 \\
    \bottomrule
    \end{tabular}%

    }
  \label{tab:addlabel}%
\end{table}%

\subsection{Power System Capacity Mix}
\begin{figure}[H]
    \centering
    \includegraphics[width = 0.8\textwidth]{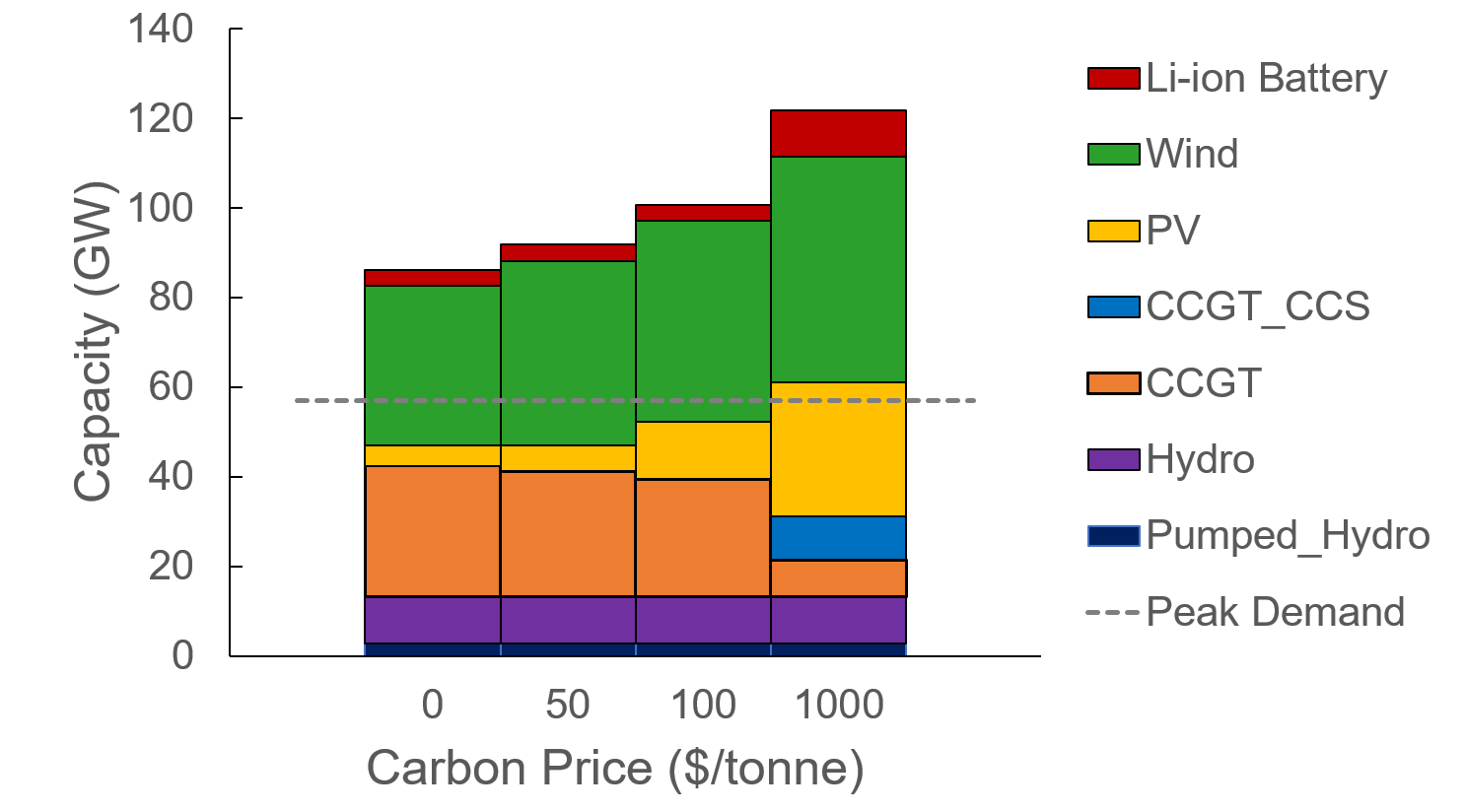}
    \caption{Optimal Power System Planning Results for the U.S. Northeast in 2050}
    \label{fig:power_system_mix}
\end{figure}

% \clearpage
\section{Hydrogen Generation Cost and Performance Parameters}

\begin{table}[!h]
  \begin{center}
  \caption{Parameters of H$_2$ Generation Technologies, in 2019 US\$ \cite{IEA2019,IEAGHG,nrel_el,NRELH2AAnalysis,PROOST20194406}}
  \setlength{\tabcolsep}{4pt}
    \begin{tabular}{lccc}
    \toprule
    \multicolumn{1}{l}{} & \multicolumn{1}{c}{Electrolysis} & \multicolumn{1}{c}{SMR} & \multicolumn{1}{c}{SMR w/ CCS} \\
    \midrule
    $\text{M}^{\text{GEN}}_{k,z} $ (tonne/hour) & 0.06   & 9.2   & 9.2 \\
    $\text{c}^{\text{GEN}}_{k} $ (M\$) & 3     & 161   & 296 \\
    $\eta^{\text{ELE}}_{k,z} $ (MWh/tonne) & 53    & 0     & 0 \\
    $\eta^{\text{GAS}}_{k,z} $ (MMBtu/tonne) & 0     & 146   & 160 \\
    $ \text{e}^{\text{GEN}}_{k} $ (tCO$_2$/tH$_2$) & 0     & 10.0  & 1.0 \\
    Lifetime (year) & 10    & 25    & 25 \\
    \bottomrule
    \end{tabular}%
  \label{tab:gen_h2}
    \end{center}
\end{table}

\clearpage
\section{Hydrogen Transmission and Storage Parameters}
\begin{table}[h]
  \begin{center}
  \caption{Parameters of H$_2$ Transmission and Storage in 2019 US\$ \cite{yang_determining_2007,schoenung11,FEKETE201510547,ingaa,samsatli_multi-objective_2018}}
  \label{tab:transmission_h2}
  \setlength{\tabcolsep}{8pt}
\begin{tabular}{lclc}
\toprule
Pipeline &       & Gas Tank &  \\
\midrule
$ \text{c}^{\text{PIP}}_{i} $ (M\$/mile) & 2.8   & $ \text{c}^{\text{STO}}_{s} $ (M\$/tonne) & 0.58 \\
$ \overline{\text{E}}_{i}^{\text{PIP}} $ (tonne/mile) & 0.3   & Lifetime (year) & 12 \\
Lifetime (year) & 40    & $ \Phi^{\text{c}}_{s} $ (\$/(tonne/hour)) & 0.5 \\
$ \alpha^{\text{c}}_{i} $ (\$/mile-unit) & 700 & $ \Phi^{\text{e}}_{s} $ (MWh/tonne)  & 2 \\ $ \beta^{\text{c}}_{i} $ (\$/unit) & 0.75  
  & \multicolumn{1}{c}{} &  \\
$ \alpha^{\text{e}}_{i} $ (MWh/tonne-mile)& 1     & \multicolumn{1}{c}{} &  \\
$ \beta^{\text{e}}_{i} $ (MWh/tonne) & 1     & \multicolumn{1}{c}{} &  \\
$ \underline{\text{R}}^{\text{PIP}}_{i} $ & 0     & \multicolumn{1}{c}{} &  \\
\midrule
Liquid Truck &       & Gas Truck &  \\
\midrule
$ \overline{\text{E}}_{j}^{\text{TRU}} $ (tonne) & 4     & $ \overline{\text{E}}_{j}^{\text{TRU}} $ (tonne) & 0.3 \\
$ \text{c}^{\text{TRU}}_{j} $ (M\$/mile) & 0.8   & $ \text{c}^{\text{TRU}}_{j} $ (M\$/mile) & 0.3 \\
$ \text{o}^{\text{TRU}}_{j} $ (\$/mile) & 1.5   & $ \text{o}^{\text{TRU}}_{j} $ (\$/mile) & 1.5 \\
$ \sigma_{j} $ & 0  & $ \sigma_{j} $  & 3\% \\
Lifetime (year) & 12    & Lifetime (year) & 12 \\
$ \Gamma^{\text{c}}_{j} $ (\$/(tonne/hour)) & 32    & $ \Gamma^{\text{c}}_{j} $ (\$/(tonne/hour)) & 1.5 \\
$ \Gamma^{\text{e}}_{j} $ (MWh/tonne) & 11    & $ \Gamma^{\text{e}}_{j} $ (MWh/tonne) & 1 \\
\bottomrule
\end{tabular}%
    \end{center}
\end{table}

% \subfile{Sections/nomenclature}
\clearpage
\bibliography{MyLibrary}
\bibliographystyle{IEEEtran}

%% file: Sections/0_nomenclature.tex
\section*{Nomenclature}
\begin{nomenclature}
    \subsection*{Indices and Sets}
    \entry{$i$, $\mathbb{I}$}{Index and set of H$ _2 $ pipelines options in the network}
    \entry{$j$, $\mathbb{J}$}{Index and set of H$_2$ truck resources}
    \entry{$k$, $\mathbb{K}$}{Index and set of H$_2$ generation resources}
    \entry{$s$, $\mathbb{S}$}{Index and set of H$_2$ storage resources}
    \entry{$t,\tau$}{Index of time intervals}
    % \entry{$i$}{Index and set of H$ _2 $ pipelines options in the network}
    % \entry{$j$}{Index of H$_2$ truck resources}
    % \entry{$k$}{Index of H$_2$ generation resources}
    % \entry{$s$}{Index of H$_2$ storage resources}
    % \entry{$t,\tau$}{Index of time intervals}
    \entry{$z,z^{'}$}{Index of zones in the network}
    \entry{$z \rightarrow z^{'}$}{Index of paths for H$ _2 $ transport}
    \entry{$\mathbb{B}$}{Set of paths for H$ _2 $ transport}
    % \entry{$\mathbb{I}$}{Set of H$ _2 $ pipelines options in the network}
    % \entry{$\mathbb{J}$}{Set of H$_2$ truck resources}
    % \entry{$\mathbb{K}$}{Set of H$_2$ generation resources}
    % \entry{$\mathbb{S}$}{Set of H$_2$ storage resources}
    \entry{$\mathbb{T}$}{Set of time intervals}
    \entry{$\mathbb{Z}$}{Set of zones in the network}
    \entry{$\text{E}$}{Superscript for empty truck}
    \entry{$\text{F}$}{Superscript for full truck}
    \entry{$\text{CHA}$}{Superscript for charging}
    \entry{$\text{DIS}$}{Superscript for discharging}
    \entry{$\text{c}$}{Superscript for capital cost}
    \entry{$\text{e}$}{Superscript for electricity consumption}
    
    \subsection*{Parameters}
    \entry{$t_{0} $}{Time origin}
    \entry{$\Delta t$}{Time resolution (hour)}
    \entry{$\delta_{()} $}{Annuity factor for various resources}
    \entry{$\text{c}^{\text{GEN}}_{k} $}{Unit capital cost of H$_{2}$ generation unit $ k $ (\$/(tonne-H$_{2}$/hour))}
    % \entry{$\lambda^{\text{ELE}}_{z,t}  $, $\lambda^{\text{GAS}}_{z,t} $}{Electricity and natural gas price at zone $ z $ and time $ t $ (\$/MWh$_{e}$, \$/MMBtu)}
    \entry{$\lambda^{\text{ELE}}_{z,t}  $}{Electricity price at zone $ z $ and time $ t $ (\$/MWh$_{e}$)}
    \entry{$\lambda^{\text{GAS}}_{z,t} $}{Natural gas price at zone $ z $ and time $ t $ (\$/MMBtu)}
    \entry{$\eta^{\text{ELE}}_{k,z} $}{Electricity consumption rate per unit H$_{2}$ production of resource $ k $ at zone $ z $ (MWh$_{e}$/tonne-H$_{2}$)}
    \entry{$\eta^{\text{GAS}}_{k,z} $}{Natural gas consumption rate per unit H$_{2}$ production of resource $ k $ at zone $ z $ (MMBtu-gas/tonne-H$_{2}$)}
    \entry{$\text{M}^{\text{GEN}}_{k,z} $}{The rated production capacity of a H$_{2}$ generation unit for resource $ k $ at zone $ z $ (tonne-H$_{2}$/hour)}
    % \entry{$\underline{\text{R}}^{\text{GEN}}_{k,z} $}{The minimum stable power input as a fraction of the rated power capacity for resource $ k $ at zone $ z $}
    % \entry{$\overline{\text{R}}^{\text{GEN}}_{k,z} $}{The maximum stable power input as a fraction of the rated power capacity for resource $ k $ at zone $ z $}
    % \entry{$\tau^\text{UP}_{k,z} $}{The minimum up time for units of resource $ k $ at zone $ z $ (hour)}
    % \entry{$\tau^\text{DOWN}_{k,z} $}{The minimum down time for units of resource $ k $ at zone $ z $ (hour)}
    \entry{$ \text{c}^{\text{STO}}_{s} $}{ Unit capital cost of H$_{2}$ storage $ s $ (\$/tonne-H$_{2}$) }
    % \entry{$ \eta^{\text{STO}}_{s,z} $}{ Charge/discharge efficiency of storage resource $ s $ at zone $ z $ (\%) }
    % \entry{$ \underline{\text{R}}^{\text{STO}}_{s,z} $}{ The minimum required H$_{2}$ fraction of storage resource $ s $ at zone $ z $ for cushion gas }
    \entry{$ \text{L}_{z \rightarrow z^{'}} $}{ Distance between zone $ z $ and $ z' $, used for calculating pipeline length and road transport distance (miles) }
    \entry{$ \overline{\text{F}}_{i} $}{ The maximum injecting/withdrawing flow rate of the pipeline $ i $ (tonne-H$_{2}$/hour) }
    % \entry{$ \overline{\text{E}}_{i}^{\text{PIP}} $}{ The capacity (the maximum volume of H$ _2 $ that can be stored) of the pipeline $ i $ per unit length (tonne-H$_{2}$/mile) }
    % \entry{$ \underline{\text{R}}^{\text{PIP}}_{i} $}{ The minimum required H$_{2}$ volume per unit of capacity of the pipeline $ i $ }
    \entry{$ \text{c}^{\text{PIP}}_{i} $}{ Unit capital cost of H$_{2}$ pipeline of type $ i $ (fixed diameter) (\$/mile) }
    \entry{$ \text{c}^{\text{TRU}}_{j} $}{ Unit capital cost of truck $ j $ (\$/unit) }
    \entry{$ \text{o}^{\text{TRU}}_{j} $}{ Unit operational cost of truck $ j $ (\$/mile)) }
    \entry{$ \overline{\text{E}}_{j}^{\text{TRU}} $}{ Capacity of truck technology $ j $ (tonne-H$_{2}$) }
    \entry{$ \sigma_{j} $}{ H$ _2 $ boil-off loss of truck technology $ j $ }
    \entry{$ \Delta_{z \rightarrow z^{'}} $}{ Travel time delay of trucks between zone $ z $ and $ z' $ }
    \entry{$ \text{e}^{\text{GEN}}_{k} $}{ Unit emission rate of H$_{2}$ generation resorce $ k $ (tonne-CO$_{2}$/tonne-H$_{2}$) }
    \entry{$ \text{e}^{\text{TRU}}_{j} $}{ Emission rate of truck $ j $ (tonne-CO$_{2}$/(tonne-H$_{2}$-mile)) }
    \entry{$ \text{c}^{\text{LOS}} $}{ Unit value of lost H$_{2}$ load (\$/(tonne/hour)) }
    \entry{$ \alpha^{\text{c}}_{i} $, $ \alpha^{\text{e}}_{i} $}{ Unit capital cost and electricity consumption of compression facilities proportional to pipeline length for pipeline type $ i $ (\$/mile, MWh$_{e}$/tonne-H$_{2}$-mile) }
    % \entry{$ \alpha^{\text{c}}_{i} $}{ Unit capital cost of compression facilities proportional to pipeline length for pipeline type $ i $ (\$/mile) }
    % \entry{$ \alpha^{\text{e}}_{i} $}{ Unit energy consumption of compression facilities proportional to pipeline length for pipeline type $ i $ (MWh$_{e}$/tonne-H$_{2}$-mile) }
    \entry{$ \beta^{\text{c}}_{i} $, $ \beta^{\text{e}}_{i} $}{ Unit capital cost and electricity consumption of compression facilities not related to pipeline length for pipeline type $ i $ (\$, MWh$_{e}$/tonne-H$_{2}$) }
    % \entry{$ \beta^{\text{c}}_{i} $}{ Unit capital cost of compression facilities not related to pipeline length for pipeline type $ i $ (\$) }
    % \entry{$ \beta^{\text{e}}_{i} $}{ Unit energy consumption of compression facilities not related to pipeline length for pipeline type $ i $ (MWh$_{e}$/tonne-H$_{2}$) }
    \entry{$ \Gamma^{\text{c}}_{j} $, $ \Gamma^{\text{e}}_{j} $}{ Unit capital cost and electricity consumption of compression/liquefaction facilities for truck type $ j $ (\$/(tonne-H$_{2}$/hour), MWh$_{e}$/tonne-H$_{2}$) }
    % \entry{$ \Gamma^{\text{c}}_{j} $}{ Unit capital cost of compression/liquefaction facilities for truck type $ j $ (\$/(tonne-H$_{2}$/hour)) }
    % \entry{$ \Gamma^{\text{e}}_{j} $}{ Unit energy consumption of compression/liquefaction facilities for truck type $ j $ (MWh$_{e}$/tonne-H$_{2}$) }
    \entry{$ \Phi^{\text{c}}_{s} $, $ \Phi^{\text{e}}_{s} $}{ Unit capital cost and electricity consumption of compression facilities for storage type $ s $ (\$/(tonne-H$_{2}$/hour), MWh$_{e}$/tonne-H$_{2}$) }
    % \entry{$ \Phi^{\text{c}}_{s} $}{ Unit capital cost of compression facilities per storage charging/discharging capacity for storage type $ s $ (\$/(tonne-H$_{2}$/hour)) }
    % \entry{$ \Phi^{\text{e}}_{s} $}{ Unit energy consumption of compression facilities for storage type $ s $ (MWh$_{e}$/tonne-H$_{2}$) }
    \entry{$ \text{D}_{z,t} $}{ H$_{2}$ demand at zone $ z $ time $ t $ (tonne-H$_{2}$/hour) }
    \entry{$ \Omega_{t} $}{ Annual scaling factor for time $ t $ }
    \entry{$ \text{c}^{\text{EMI}} $}{ Carbon emission price (\$/tonne) }
    
    \subsection*{Variables}
    \entry{$ h^{\text{GEN}}_{k,z,t} $}{ H$_2$ produced by resource $ k $ at zone $ z $ during time $ t $ (tonne/hour) }
    \entry{$ h^{\text{TRA}}_{z,t} $}{ H$_2$ transport to zone $ z $ during time $ t $ (tonne/hour) }
    \entry{$ h^{(\cdot)}_{s,z,t} $}{ H$_2$ consumed to charge or released by storage resource $ s $ at zone $ z $ during time $ t $ (tonne/hour) }
    % \entry{$ h^{(\cdot)}_{s,z,t} $}{ H$_2$ released by storage resource $ s $ at zone $ z $ during time $ t $ (tonne/hour) }
    \entry{$ h^{\text{LOS}}_{z,t} $}{ Lost H$_{2}$ load at zone $ z $ during time $ t $ (tonne/hour) }
    \entry{$ g^{\text{GAS}}_{k,z,t} $}{ Natural gas used by resource $ k $ to produce H$_{2}$ at zone $ z $ during time $ t $ (MMBtu/hour) }
    % \entry{$ n_{k,z,t} $}{ Number of online units of resource $ k $ at zone $ z $ during time $ t $ }
    % \entry{$ n^{\text{UP}}_{k,z,t} $}{ Number of starting-up units of resource $ k $ at zone $ z $ during time $ t $ }
    % \entry{$ n^{\text{DOWN}}_{k,z,t} $}{ Number of shutting-down units of resource $ k $ at zone $ z $ during time $ t $ }
    \entry{$ N_{k,z} $}{ Number of available units of resource $ k $ at zone $ z $ }
    % \entry{$ N^{\text{existing}}_{k,z} $}{ Number of existing units of resource $ k $ at zone $ z $ }
    % \entry{$ N^{\text{new}}_{k,z} $}{ Number of newly installed units of resource $ k $ at zone $ z $ }
    % \entry{$ N^{\text{retired}}_{k,z} $}{ Number of retired units of resource $ k $ at zone $ z $ }
    \entry{$ V^{\text{STO}}_{s,z} $}{ Capacity of H$_{2}$ storage resource $ s $ at zone $ z $ (tonne) }
    \entry{$ H^{\text{STO}}_{s,z} $}{ Maximum charge/discharge rate of H$_{2}$ storage resource $ s $ at zone $ z $ (tonne/hour) }
    \entry{$ l_{z \rightarrow z',i} $}{ Number of pipeline type $ i $ between $ z' $ and $ z $ }
    \entry{$ h^{\text{PIP}}_{z \rightarrow z',i,t} $}{ H$ _2 $ exchange at zone $ z $ via pipeline $ i $ between $ z' $ and $ z $ during time $ t $ (tonne/hour). $ \text{PIP+} $ is the superscript for delivered H$ _2 $, and $ \text{PIP-} $ for drawn H$ _2 $}
    % \entry{$ h^{\text{PIP+}}_{z \rightarrow z',i,t} $}{ H$ _2 $ delivered to zone $ z $ via pipeline $ i $ between $ z' $ and $ z $ during time $ t $ (tonne/hour) }
    % \entry{$ h^{\text{PIP-}}_{z \rightarrow z',i,t} $}{ H$ _2 $ drawn from zone $ z $ via pipeline $ i $ between $ z' $ and $ z $ during time $ t $ (tonne/hour) }
    \entry{$ p^{\text{COM}}_{z,t} $}{ Power consumption by compression facilities at zone $ z $ during time $ t $ (MW) }
    \entry{$ H^{\text{TRU}}_{z,j} $}{ Maximum compression/liquefaction capacity of truck station type $ j $ at zone $ z $(tonne/hour) }
    % \entry{$ h^{\text{TRU}}_{z,j,t} $}{ H$ _2 $ exchange at zone $ z $ via truck of type $ j $ during time $ t $ (tonne/hour) }
    \entry{$ v^{(\cdot)}_{j,t} $}{ Number of full or empty trucks of type $ j $ at time $ t $ }
    % \entry{$ v^{\text{E}}_{j,t} $}{ Total number of empty trucks of type $ j $ at time $ t $ }
    \entry{$ V_{j} $}{ Total number of invested trucks of type $ j $ }
    \entry{$ u^{(\cdot)}_{z \rightarrow z^{'},j,t} $}{ Number of full or empty trucks of type $ j $ in transit from $ z $ to $ z' $ at time $ t $ }
    % \entry{$ u^{\text{E}}_{z \rightarrow z^{'},j,t} $}{ Number of empty trucks of type $ j $ in transit from $ z $ to $ z' $ at time $ t $ }
    \entry{$ x^{(\cdot)}_{z \rightarrow z^{'},j,t} $}{ Number of full or empty trucks of type $ j $ departing from $ z $ to $ z' $ at time $ t $ }
    % \entry{$ x^{\text{E}}_{z \rightarrow z^{'},j,t} $}{ Number of empty trucks of type $ j $ departing from $ z $ to $ z' $ at time $ t $ }
    \entry{$ y^{(\cdot)}_{z \rightarrow z^{'},j,t} $}{ Number of full or empty trucks of type $ j $ arriving at $ z' $ from $ z $ at time $ t $ }
    % \entry{$ y^{\text{E}}_{z \rightarrow z^{'},j,t} $}{ Number of empty trucks of type $ j $ arriving at $ z' $ from $ z $ at time $ t $ }
    \entry{$ q^{(\cdot)}_{z,j,t} $}{ Number of full, empty, charged, or discharged trucks of type $ j $ available at $ z $ at time $ t $ }
    % \entry{$ q^{\text{E}}_{z,j,t} $}{ Number of empty trucks of type $ j $ available at $ z $ at time $ t $ }
    % \entry{$ q^{\text{CHA}}_{z,j,t} $}{ Number of trucks of type $ j $ that has just been fully charged at $ z $ at time $ t $ }
    % \entry{$ q^{\text{DIS}}_{z,j,t} $}{ Number of trucks of type $ j $ that has just been fully discharged at $ z $ at time $ t $ }

\end{nomenclature}

%% file: Sections/1_introduction.tex
\section{Introduction}
% \section{Literature Review}
Deep decarbonization of the energy system is contingent on identifying pathways for eliminating greenhouse gas (GHG) emissions from not only the power sector but also other end-use sectors where direct electrification may be challenging \cite{davis2018net}. In this context, identifying cost-effective pathways for supplying energy carriers like hydrogen remains an appealing prospect \cite{IEA2019}. Recent renewed interest in hydrogen has been spurred, in part, by expectations on cost declines for water electrolyzers \cite{Schmidt17}, which raises the prospect of electrolytic hydrogen produced from variable renewable energy (VRE) resources becoming cost-competitive with fossil-fuel based pathways such as steam methane reforming (SMR) \cite{GUERRA20192425}. However, hydrogen production represents only a fraction of the total cost of hydrogen supply for distributed end uses like transportation, owing to the relatively high cost associated with transmission, storage, and distribution \cite{DEMIR201810420}. Therefore, identification of cost-effective hydrogen supply chains (HSC) requires a careful consideration of all stages of the supply chain, including production, transport, storage and end-use, as well as their inter-dependencies.

Some existing studies have investigated hydrogen as a basic storage system (with co-located electrolyzer, fuel cell, and gas tank) while not as a network that consists of various production, transmission and storage infrastructures as well as hydrogen demand. The coordination schemes of hydrogen storage with wind farms are studied in \cite{alkano18,recalde14}. The scheduling and energy trading strategies of hydrogen storage in wholesale or local energy markets are studied in \cite{xiao18,el-taweel19}.
% On top of the basic hydrogen storage system setting, the power to heat and hydrogen technology is modelled in \cite{fu19,pan20,yun19}, and the hydrogen enriched compressed natural gas is incorporated in \cite{pan20}, to consider the interactions among power, heat, and gas sectors. 
% Hydrogen storage is compared with other types of storage such as battery storage, pumped hydrogen storage, thermal energy storage, etc. as options to decarbonize European energy system in \cite{victoria19}. 
These studies revealed the potential role and benefit of hydrogen storage in power systems, while they did not explore the economic competitiveness of power-to-hydrogen technologies nor inform the investment needed in HSC infrastructures for hydrogen use in power and other end-use sectors.
 
In the literature on HSC planning, many opportunities for further modeling improvements exist. First, few HSC models account for the variability in electricity supply and its impact on the least-cost infrastructure planning. With the prospect of electrolytic hydrogen supply, there is a need to incorporate the arising spatiotemporal variations in electricity supply (and prices) resulting from increasing VRE penetrations in the power systems \cite{KETTERER2014270,BRANCUCCIMARTINEZANIDO2016474}. Second, most existing literature on HSC planning adopt a limited set of hydrogen transmission pathways, ignoring pipeline \cite{li_19,almansoori16,kim17, ochoa18,de-leon_almaraz_hydrogen_2014}, gas truck \cite{de-leon_almaraz_hydrogen_2014}, and/or liquid truck \cite{kim17} options. Some also ignore critical hydrogen production options including electrolyzer \cite{almansoori16} and SMR with or without carbon capture and storage (CCS) \cite{emonts_19,kim17,li_19,reu17,welder_18}, even though the mode of hydrogen production (electrolyzer versus SMR) is a critical determining factor to the cost-competitive transmission mode (gas/liquid truck versus pipeline). Third, most literature use over-simplified models for hydrogen transmission via truck transport. For example, the truck transportation is modelled with fixed lower and upper flow limits for each route in \cite{almansoori16,kim17,li20,ochoa18,reu17,welder_18}. In reality, traveling delay is inevitable for road transport, and a truck can serve different routes and sites both as transmission and storage assets. The availability of trucks is also likely to vary over space and time, and assuming the hydrogen flow limits are time- and route-invariant will underestimate either the required capital costs or the flexibility of trucks. To the best of our knowledge, there has not been a full-fledged HSC model that incorporates all critical technological options in hydrogen production, transmission, and storage while also accounting for their spatiotemporal operational flexibility. 
  
This paper develops a high-fidelity HSC model to evaluate the cost-competitiveness of hydrogen for decarbonizing various end-uses and the trade-offs between various technology options across the HSC. The main contributions of this paper are summarized as follows:

1) We propose an HSC planning model that determines the least-cost technology mix across the supply chain based on its operation under spatiotemporal variations in electricity inputs and hydrogen demands. 

2) We explicitly model the flexibility of a wide range of conventional and emerging hydrogen-related technologies, including electrolyzer, SMR with and without CCS, pressure vessel and geological storage, hydrogen transport via liquid and gaseous modes in trucks, gaseous mode via pipeline, and associated compression/liquefaction facilities. 

3) The framework incorporates a flexible scheduling and routing model for hydrogen trucks to serve as both transmission and mobile storage, which shift hydrogen demand/production in space and time while being shared across the whole hydrogen network. We also develop a scheduling model for pipeline operation that considers the ability of pipeline to function as storage via an approximation of line packing. While the HSC planning model with flexible truck scheduling is complex in nature, we made the model computationally tractable using judicious approximations at little cost of modelling accuracy. 

4) We apply the HSC planning model to study the future hydrogen infrastructure needs under various carbon policy and hydrogen demand scenarios in the U.S. Northeast, based on electricity price profiles generated from a power system capacity expansion model. We demonstrate the value of the proposed sophisticated scheduling models and analyze the trade-offs and synergies between different technologies given their investment and operational flexibility.

The remainder of this paper is organized as follows. Section \ref{flexibility_modelling} introduces the flexibility characteristics of various transmission and storage technologies across the HSC. Section \ref{model_formulation} formulates the HSC planning model in detail. Section IV describes the case studies and the results of the HSC planning model for various technology and policy scenarios. The conclusions are drawn in Section \ref{conclusion}.

%% file: Sections/2_flexibility_modelling.tex
\section{Flexibility of Storage and Transmission Technologies in Hydrogen Supply Chain}
\label{flexibility_modelling}
H$_2$ storage resources, both stationary and mobile storage either in geological formations and pressure vessels, enable shifting production in space and time to better match supply and demand. Among stationary storage, we focus on physical H$_2$ storage technologies, namely pressure vessel H$_2$ storage and geological H$_2$ storage, which are the most mature technologies to store H$_2$ today \cite{Barthelemy17}.

% \subsection{Underground storage}
 
Geological H$_2$ storage is based on storing thousands of tonnes of compressed H$_2$ in geological formations like salt caverns and aquifers. The limited geographic availability of these formations and perceived safety issues, however, may limit their use. Pressure vessel H$_2$ storage, unlike geological storage, does not depend on geographical resources and thus can be deployed in a modular and spatially flexible manner. The cost of pressure vessel H$_2$ storage at pressures of 350-500 bar is low compared to electrochemical energy storage, typically ranging from \$15-20/kWh-H$_2$ \cite{yang_determining_2007,schoenung11,Barthelemy17}.
% and is currently practiced at 5 locations globally at 1000's of tonnes of storage capacity \cite{drive2017hydrogen}. 
Cushion gas is required to maintain the pressure and withdraw rate of geologic storage, typically ranging from 30-50\% of the total storage volume, depending on the geological formation type \cite{Lord14}.
% Estimated capital costs of geological H$_2$ storage are less than \$1/kWh-H$_2$ \cite{schoenung11,Lord14}. However, the limited geographic availability of these sites and perceived safety issues may limit their use.  
% \subsection{Aboveground gas tank}

% \subsection{Gas truck}
Among mobile H$_2$ storage options, compressed H$_2$ trucks, or namely gas trucks, carry pressure vessels (bundles of H$_2$ tubes) with hydrogen usually compressed at 180 bar, much lower than the pressure vessel stationary storage (and hence lower carrying weight). Due to the pressure limit imposed by road transport regulations, the hydrogen carrying capacity of a gas truck is limited to 0.3-0.6 tonne (10-17 MWh$_\text{t}$) \cite{drive2013hydrogen,yang_determining_2007}. The advantage of gas trucks over pressure vessel H$_2$ storage is that it can travel and be shared across the whole HSC to provide on-demand storage service and meet spatio-temporally varying energy supply and demand.
% \subsection{Liquid truck}
Liquid trucks carry cryogenic vessels with liquid hydrogen (lower than 33 K). The hydrogen carrying capacity of commercially available liquid truck is typically 4 tonnes (133 MWh$_\text{t}$) \cite{drive2013hydrogen,yang_determining_2007}, significantly greater than gas trucks owing to higher volumetric energy density of liquid hydrogen. Liquefaction, however, requires much higher energy consumption and capital costs than gas compression (approximately 11 MWh/tonne vs. 1 MWh/tonne \cite{yang_determining_2007}). Liquid trucks also have higher boil-off rates and thus are not a preferred solution for long-duration storage but may be more suitable for long-distance transport. 

% \subsection{Pipeline}
While pipelines are usually built for large-volume transmission, it can also function as storage, through line packing. Hydrogen can be either withdrawn from or injected to the pipeline terminals like storage discharging or charging. A pipeline with 8-inch diameter and 100 bar can store approximately 0.3 tonne hydrogen per mile, approximately 20\% of the energy content of natural gas that can be stored in a similar natural gas pipeline per mile \cite{oseghale2011application}.  

% Specialized materials of constructions are needed when handling pressurized H$_2$ due to H$_2$ embrittlement and related fatigue. This could both increase cost relative to cost of handling pressurized natural gas and limit deployment due to perceived safety issues. For example, our reference estimate of H$_2$ pipeline costs are approximately 5 times the cost of a comparable natural gas pipeline after accounting for volumetric energy density differences between H$_2$ and natural gas. 

% The two commercially available water electrolyzer technologies, proton exchange membrane and alkaline electrolyzer, are both capable of adjusting their electricity inputs rapidly and operating over a wide range of utilization levels (25\%-100\% rated capacity) which makes them suitable to integrate with VRE sources \cite{Eichman16,Glenk19}. Their sizes are highly configurable due to their modular design, from several kWs to MWs \cite{Ursua12}, and thus can be deployed in both distributed and centralized manners.

% Hydrogen production technologies of natural gas reforming, such as SMR and autothermal reforming (ATR), are usually deployed as large centralized plants to achieve economies of scale with limited operational flexibility. With an operating temperature of $900-1000{\celsius}$, these natural gas reforming technologies require a minimum start time to reach the temperature and a minimum operating level to maintain the temperature \cite{catal9100801}.

%% file: Sections/3_model_formulation.tex
\section{Hydrogen Supply Chain Planning Model}
\label{model_formulation}
\begin{figure}[t]
		\centering
		\includegraphics[width=0.48\textwidth]{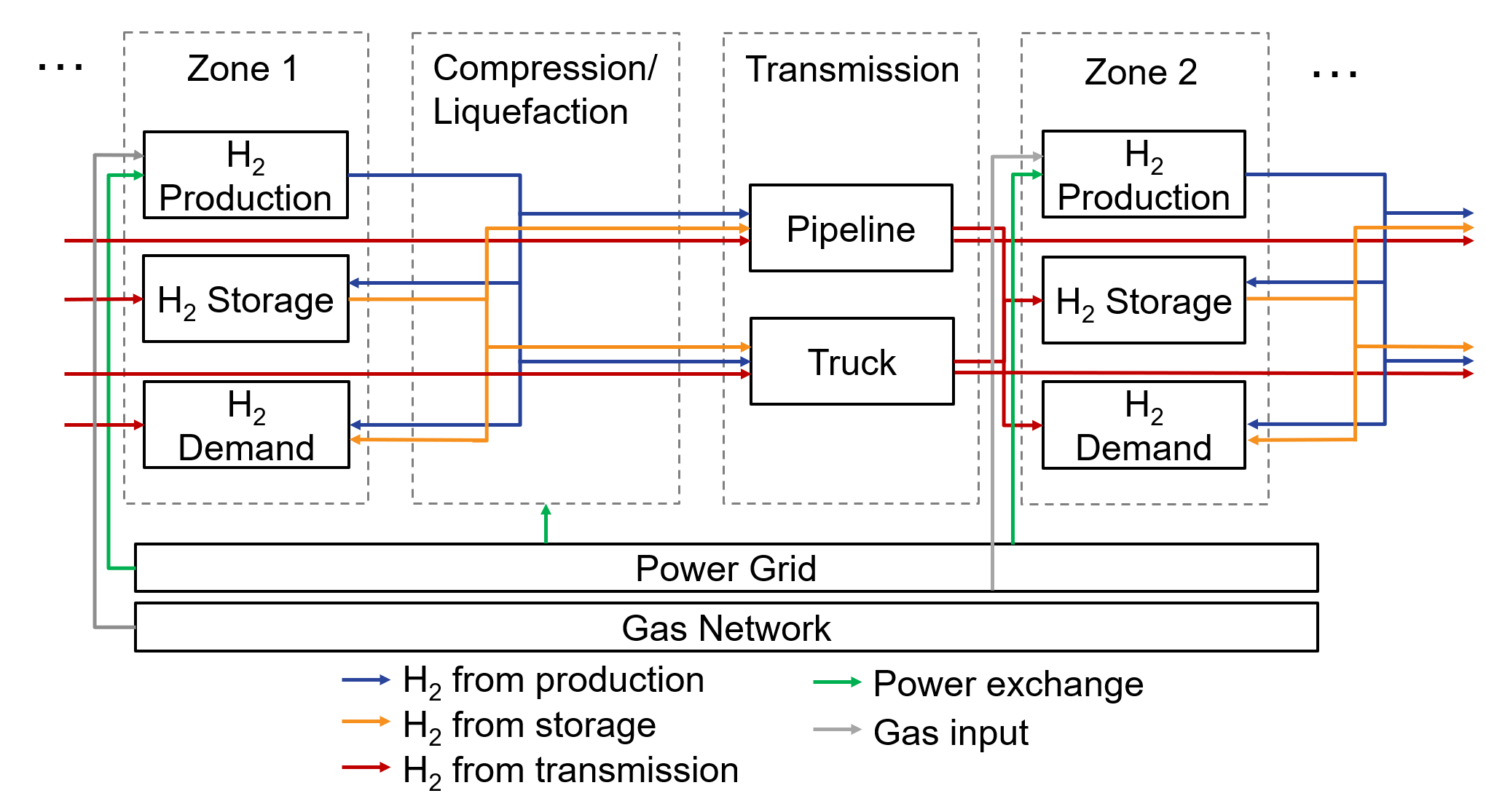}
		\caption{Schematic of the hydrogen supply chain model.}
		\label{fig-schematic}
\end{figure}

Here we present the model framework for identifying the least-cost H$_2$ infrastructure, spanning production, transmission and storage, needed to meet given spatio-temporal H$_2$ demand and energy (electricity, natural gas) prices.  The schematic of the HSC including various components and energy flows are shown in Figure \ref{fig-schematic}. The interactions between the HSC and the power grid are modeled assuming that the HSC is a price-taker, which implies that electricity price profiles at a given location are not impacted by electricity consumption from the HSC at that location. The decision variables are the capacities/numbers of units for H$_2$ production, storage, compression, and transmission resources, as well as the hourly operational schedules of the resources. They are continuous except the truck scheduling variables.
% Investments in SMR with and without CCS are modeled as integer variables, to consider the effect of economies of scale and their relative flexibility limitations. All other decision variables are continuous, which constitute a mixed integer linear programming model.

\subsection{Objective function}
The total cost associated with H$_2$ infrastructure includes the capital and operating cost of the four main elements of the supply chain, production, storage, compression and transmission (as in (\ref{Total_cost})). This includes the capital costs of: H$_2$ production units ($C_{\text{GEN}}^{\text{c}}$) as in (\ref{Production_cost}), H$_2$ storage ($C_{\text{STO}}^{\text{c}}$) as in (\ref{Storage_cost}), H$_2$ pipeline ($C_{\text{PIP}}^{\text{c}}$) as in (\ref{Pipeline_cost}), trucks for H$_2$ transport ($C_{\text{TRU}}^{\text{c}}$) as in (\ref{Truck_cost}), and H$_2$ compression facilities ($C_{\text{COM}}^{\text{c}}$) associated with pipeline ($C_{\text{COM,PIP}}^{\text{c}}$), truck station ($C_{\text{COM,TRU}}^{\text{c}}$), and storage ($C_{\text{COM,STO}}^{\text{c}}$)  as in (\ref{Compression_cost}). The main operating costs included in the objective function are: electricity cost ($C_{\text{ELE}}^{\text{o}}$) as in (\ref{Electricity_cost}), fuel i.e. natural gas, cost for H$_2$ production ($C_{\text{GAS}}^{\text{o}}$) as in (\ref{Gas_cost}),
operating costs for H$_2$ truck transport ($C_{\text{TRU}}^{\text{o}}$) as in (\ref{Truck_opex}), emission costs from H$_2$ production and truck transport ($C_{\text{EMI}}^{\text{o}}$) as in (\ref{Emission_opex}), and costs of lost load ($C_{\text{LOS}}^{\text{o}}$) as in (\ref{VOLL_opex}). 
\begin{equation}
\begin{split}
\label{Total_cost}
\min \qquad & C_{\text{GEN}}^{\text{c}} 
% + C_{\text{POW}}^{\text{c}} 
+ 
C_{\text{STO}}^{\text{c}} 
+ C_{\text{PIP}}^{\text{c}}   
+  C_{\text{TRU}}^{\text{c}} 
+ C_{\text{COM}}^{\text{c}} 
\\ + & C_{\text{ELE}}^{\text{o}} 
+ C_{\text{GAS}}^{\text{o}} 
+ C_{\text{TRU}}^{\text{o}} 
% + C_{\text{COM}}^{\text{o}} 
+ C_{\text{EMI}}^{\text{o}} 
+ C_{\text{LOS}}^{\text{o}}
\end{split}
\end{equation}

\begin{equation}
\begin{split}
\label{Production_cost}
C_{\text{GEN}}^{\text{c}} = \delta^{\text{GEN}}_{k}  \sum\limits_{z \in \mathbb{Z}}\sum\limits_{k \in \mathbb{K} }  \text{c}^{\text{GEN}}_{k}  \text{M}^{\text{GEN}}_{k,z}  N^{\text{new}}_{k,z}
\end{split}
\end{equation}

% \begin{equation}
% \begin{split}
% \label{Power_gen_cost}
% C_{\text{POW}}^{\text{c}} = \delta_{\text{POW}}  \sum\limits_{z \in \mathbb{Z}}\sum\limits_{g \in \mathbb{G}} vM^{new}_{g,z}  \text{c}^{\text{POW}}_{k}
% \end{split}
% \end{equation}

\begin{equation}
\begin{split}
\label{Storage_cost}
C_{\text{STO}}^{\text{c}} = \delta^{\text{STO}}_{s}  \sum\limits_{z \in \mathbb{Z}}\sum\limits_{s \in \mathbb{S} }  \text{c}^{\text{STO}}_{s} V^{\text{STO}}_{s,z}
\end{split}
\end{equation}

\begin{equation}
\begin{split}
\label{Pipeline_cost}
C_{\text{PIP}}^{\text{c}}  =& \delta^{\text{PIP}}_{i} \sum\limits_{z \rightarrow z^{'} \in \mathbb{B}} \sum\limits_{i \in \mathbb{I}} \text{c}^{\text{PIP}}_{i}   \text{L}_{z \rightarrow z^{'}}     l_{z \rightarrow z',i}
\end{split}
\end{equation}

\begin{equation}
\begin{split}
\label{Truck_cost}
C_{\text{TRU}}^{\text{c}} = \delta^{\text{TRU}}_{j}  \sum\limits_{j \in \mathbb{J}} \text{c}^{\text{TRU}}_{j} V_{j}  
\end{split}
\end{equation}

\begin{equation}
\label{Compression_cost}
    \begin{split}
        C_{\text{COM}}^{\text{c}} =& \delta_{\text{COM}}  
        ( C_{\text{COM,PIP}}^{\text{c}} + C_{\text{COM,TRU}}^{\text{c}} + C_{\text{COM,STO}}^{\text{c}} )
        \\=& \delta_{\text{COM}} \left[ \sum\limits_{z \rightarrow z^{'} \in \mathbb{B}} \sum\limits_{i \in \mathbb{I}}  \left(  \alpha^{\text{c}}_{i}  \text{L}_{z \rightarrow z'} + \beta^{\text{c}}_{i} \right) l_{z \rightarrow z',i} \right.   
        \\+& \left.\sum\limits_{z \in \mathbb{Z}} \sum\limits_{j \in \mathbb{J}} \Gamma^{\text{c}}_{j} H^{\text{TRU}}_{z,j} 
        + \sum\limits_{z \in \mathbb{Z}} \sum\limits_{s \in \mathbb{S}} \Phi^{\text{c}}_{s} H^{\text{STO}}_{s,z}  \right]
    \end{split}
\end{equation}

% \begin{equation}
% \begin{split}
% \label{}
% C_{\text{COM,PIP}}^{\text{c}} =& \sum\limits_{z \rightarrow z^{'} \in \mathbb{B}} \sum\limits_{i \in \mathbb{I}} l_{z \rightarrow z',i}  \left(  \alpha^{\text{c}}_{i}  \text{L}_{z \rightarrow z'} + \beta^{\text{c}}_{i} \right)    
% \end{split}
% \end{equation}

% \begin{equation}
% \begin{split}
% \label{}
% C_{\text{COM,TRU}}^{\text{c}} =& 
% \sum\limits_{z \in \mathbb{Z}} \sum\limits_{j \in \mathbb{J}} H^{\text{TRU}}_{z,j} \Gamma^{\text{c}}_{j}
% \end{split}
% \end{equation}

% \begin{equation}
% \begin{split}
% \label{}
% C_{\text{COM,STO}}^{\text{c}} =& 
% \sum\limits_{z \in \mathbb{Z}} \sum\limits_{s \in \mathbb{S}} H^{\text{STO}}_{s,z} \Phi^{\text{c}}_{s}  
% \end{split}
% \end{equation}

\begin{equation}
\begin{split}
\label{Electricity_cost}
C_{\text{ELE}}^{\text{o}} = \sum\limits_{z \in \mathbb{Z}}\sum\limits_{t \in \mathbb{T}} \left[ \Omega_{t} \lambda^{\text{ELE}}_{z,t}  \left( \sum\limits_{k \in \mathbb{K} }   \eta^{\text{ELE}}_{k,z}  h^{\text{GEN}}_{k,z,t}
% -  \sum\limits_{k \in \mathbb{G} } vPfG_{g,z,t}
+ p^{\text{COM}}_{z,t} \right) \right] 
\end{split}
\end{equation}

\begin{equation}
\begin{split}
\label{}
p^{\text{COM}}_{z,t} = &\sum\limits_{z^{'} \in \mathbb{Z}} \sum\limits_{i \in \mathbb{I}} \left(  \alpha^{\text{e}}_{i}  \text{L}_{z \rightarrow z'} + \beta^{\text{e}}_{i} \right)    
( h^{\text{PIP+}}_{z \rightarrow z',i,t} + h^{\text{PIP-}}_{z \rightarrow z',i,t} ) 
\\+&
\sum\limits_{z \in \mathbb{Z}} \sum\limits_{j \in \mathbb{J}} \sum\limits_{t \in \mathbb{T}} \Gamma^{\text{e}}_{j}  q^{\text{CHA}}_{z,j,t}  \overline{\text{E}}_{j}^{\text{TRU}}
+
\sum\limits_{z \in \mathbb{Z}} \sum\limits_{s \in \mathbb{S}} \sum\limits_{t \in \mathbb{T}} \Phi^{\text{e}}_{s}  h^{\text{CHA}}_{s,z,t}
\end{split}
\end{equation}

\begin{equation}
\begin{split}
\label{Gas_cost}
C_{\text{GAS}}^{\text{o}} = \sum\limits_{z \in \mathbb{Z}}\sum\limits_{k \in \mathbb{K} } \sum\limits_{t \in \mathbb{T}}  \Omega_{t} \lambda^{\text{GAS}}_{z,t} h^{\text{GEN}}_{k,z,t}  \eta^{\text{GAS}}_{k,z} 
\end{split}
\end{equation}

\begin{equation}
\begin{split}
\label{Truck_opex}
C_{\text{TRU}}^{\text{o}} = &  \sum\limits_{z \rightarrow z^{'} \in \mathbb{B}} \sum\limits_{j \in \mathbb{J}} \sum\limits_{t \in \mathbb{T}} 
\Omega_{t}
\text{o}^{\text{TRU}}_{j} \text{L}_{z \rightarrow z'} \left( y^{\text{F}}_{z \rightarrow z',j,t} + 
 y^{\text{E}}_{z \rightarrow z',j,t} \right)  
\end{split}
\end{equation}

% \begin{equation}
% \begin{split}
% \label{OpexCostCompression}
% C_{\text{COM}}^{\text{o}} = &\sum\limits_{z \rightarrow z^{'} \in \mathbb{B}} \sum\limits_{i \in \mathbb{I}} H2CompressPipeUnitVarOpex_{i}  \text{L}_{z \rightarrow z^{'}}
% \\ & ( h^{\text{PIP+}}_{z' \rightarrow z,i,t} + h^{\text{PIP-}}_{z' \rightarrow z,i,t} 
% \\+& h^{\text{PIP+}}_{z \rightarrow z',i,t} + h^{\text{PIP-}}_{z \rightarrow z',i,t} ) 
% \\+&
% \sum\limits_{z \rightarrow z^{'} \in \mathbb{B}} \sum\limits_{i \in \mathbb{I}} \sum\limits_{t \in \mathbb{T}} H2CompressPipeFrontVarOpex_{i} 
% \\ & (h^{\text{PIP-}}_{z' \rightarrow z,i,t} + h^{\text{PIP-}}_{z \rightarrow z',i,t})
% \\+&
% \sum\limits_{z \in \mathbb{Z}} \sum\limits_{s \in \mathbb{S}} \sum\limits_{t \in \mathbb{T}} H2CompressTruckVarOpex_{k}  q^{\text{CHA}}_{z,j,t}  \overline{\text{E}}_{j}^{\text{TRU}}
% \\+&
% \sum\limits_{z \in \mathbb{Z}} \sum\limits_{s \in \mathbb{S}} \sum\limits_{t \in \mathbb{T}} H2CompressStorVarOpex_{k}  h^{\text{CHA}}_{s,z,t}
% \end{split}
% \end{equation}

\begin{equation}
\begin{split}
\label{Emission_opex}
C_{\text{EMI}}^{\text{o}} =&      \sum\limits_{k \in \mathbb{K}} \sum\limits_{z \in \mathbb{Z}} \sum\limits_{t \in \mathbb{T}} \Omega_{t} \text{c}^{\text{EMI}}
\text{e}^{\text{GEN}}_{k} h^{\text{GEN}}_{k,z,t}
\\+&
\sum\limits_{z \rightarrow z^{'} \in \mathbb{B}} \sum\limits_{j \in \mathbb{J}} \sum\limits_{t \in \mathbb{T}} 
\Omega_{t} \text{c}^{\text{EMI}} 
\text{e}^{\text{TRU}}_{j} \overline{\text{E}}_{j}^{\text{TRU}} \text{L}_{z \rightarrow z'} \left( y^{\text{F}}_{z \rightarrow z',j,t} + 
 y^{\text{E}}_{z \rightarrow z',j,t} \right)  
\end{split}
\end{equation}

\begin{equation}
\begin{split}
\label{VOLL_opex}
C_{\text{LOS}}^{\text{o}} = \sum\limits_{z \in \mathbb{Z}} \sum\limits_{t \in \mathbb{T}} \Omega_{t} \text{c}^{\text{LOS}} h^{\text{LOS}}_{z,t}
\end{split}
\end{equation}

% \subsection{Model Constraints}

\subsection{Hydrogen Balance Constraints}
For each zone $ z $ and time $ t $, the sum of the amount of produced H$_2$, the amount of delivered H$_2$ (net imports), and the amount of H$_2$ discharged from storage should be equal to the sum of the amount of H$_2$ charged to storage, and the H$_2$ demand minus the lost demand.
\begin{equation}
\begin{split}
\sum\limits_{k \in \mathbb{K}} h^{\text{GEN}}_{k,z,t} 
+ 
h^{\text{TRA}}_{z,t}
+ 
\sum\limits_{s \in \mathbb{S}} h^{\text{DIS}}_{s,z,t} = 
% \\&
%\sum\limits_{z' \in N,z \in \mathbb{Z}}\sum\limits_{m \in M} vH2T_{z \rightarrow z',t}
%+ 
\sum\limits_{s \in \mathbb{S}} h^{\text{CHA}}_{s,z,t}
% + 
% \sum\limits_{g \in G_{H2}} vH2P_{g,z,t}
&+ 
\text{D}_{z,t} - h^{\text{LOS}}_{z,t}
\\\qquad 
&\forall z \in \mathbb{Z}, t \in \mathbb{T}
\end{split}
\end{equation}

\subsection{Hydrogen Production Constraints}
% The amounts of electricity and gas used for hydrogen production, $ vP2G_{k,z,t} $ and $g^{\text{GAS}}_{k,z,t}$, are equal to the amount of hydrogen production, $ h^{\text{GEN}}_{k,z,t} $, times the energy consumption rate and fuel consumtion rate, $ \eta^{\text{ELE}}_{k,z} $ and $\eta^{\text{GAS}}_{k,z}$, respectively.

% \begin{equation}
% \begin{split}
% vP2G_{k,z,t} = h^{\text{GEN}}_{k,z,t}  \eta^{\text{ELE}}_{k,z}
% \qquad \forall  k \in \mathbb{K},z \in \mathbb{Z}, t \in \mathbb{T}
% \end{split}
% \end{equation}
% \begin{equation}
% \begin{split}
% g^{\text{GAS}}_{k,z,t} = h^{\text{GEN}}_{k,z,t}  \eta^{\text{GAS}}_{k,z} 
% \qquad \forall  k \in \mathbb{K},z \in \mathbb{Z}, t \in \mathbb{T}
% \end{split}
% \end{equation}

% \subsubsection{Minimum and Maximum Outputs}
The outputs of each type of H$_2$ generation facilities have to be kept within their lower and upper bounds ($ \underline{\text{R}}^{GEN}_{k,z} $ and $ \overline{\text{R}}^{GEN}_{k,z} $), as in (\ref{H2gen_lb}). $ \text{M}^{\text{GEN}}_{k,z} $ is the rated size of a H$_2$ generation unit. $ n_{k,z,t} $ denotes the number of online units. The number of online units have to be less than the available number of generation units as in (\ref{H2gen_online}). 
% The available number of units is equal to the existing units plus the new constructed units minus the retired units, as in (\ref{H2gen_number}).
\begin{equation}
\begin{split}
\label{H2gen_lb}
\overline{\text{R}}^{\text{GEN}}_{k,z}  \text{M}^{\text{GEN}}_{k,z}  n_{k,z,t} \geq h^{\text{GEN}}_{k,z,t} \geq \underline{\text{R}}^{\text{GEN}}_{k,z}  \text{M}^{\text{GEN}}_{k,z}  n_{k,z,t}
\\\qquad \forall  k \in \mathbb{K},z \in \mathbb{Z}, t \in \mathbb{T}
\end{split}
\end{equation}
% \begin{equation}
% \begin{split}
% \label{H2gen_ub}
% h^{\text{GEN}}_{k,z,t} \leq \overline{\text{R}}^{{\text{GEN}}}_{k,z}  \text{M}^{\text{GEN}}_{k,z}  n_{k,z,t}
% \qquad \forall  k \in \mathbb{K},z \in \mathbb{Z}, t \in \mathbb{T}
% \end{split}
% \end{equation}
\begin{equation}
\begin{split}
\label{H2gen_online}
n_{k,z,t} \leq N_{k,z}
\qquad \forall  k \in \mathbb{K},z \in \mathbb{Z}, t \in \mathbb{T}
\end{split}
\end{equation}
% \begin{equation}
% \begin{split}
% \label{H2gen_number}
% N_{k,z} = N^{\text{existing}}_{k,z} + N^{\text{new}}_{k,z} -  N^{\text{retired}}_{k,z}
% \qquad \forall  k \in \mathbb{K},z \in \mathbb{Z}
% \end{split}
% \end{equation}

% \subsubsection{Minimum Up and Down Times}

The number of units starting up and shutting down, represented by $n^{\text{UP}}_{k,z,t}$ and $n^{\text{DOWN}}_{k,z,t}$, are modelled in (\ref{H2gen_start_shut}). There are limits on the period of time between when a unit starts up and when it can be shut-down again, and vice versa, modelled in (\ref{H2gen_start_up}) and (\ref{H2gen_shut_down}). The minimum up and down time are denoted by $ \tau^\text{UP}_{k,z} $ and $ \tau^\text{DOWN}_{k,z} $, respectively.
\begin{equation}
\begin{split}
\label{H2gen_start_shut}
n_{k,z,t} - n_{k,z,t-1} = n^{\text{UP}}_{k,z,t} - n^{\text{DOWN}}_{k,z,t}
\qquad \forall  k \in \mathbb{K},z \in \mathbb{Z}, t \in \mathbb{T}
\end{split}
\end{equation}
\begin{equation}
\begin{split}
\label{H2gen_start_up}
n_{k,z,t} \geq 
\sum\limits_{\tau = t-\tau^\text{UP}_{k,z}}^{t} n^{\text{UP}}_{k,z,t}
\qquad \forall  k \in \mathbb{K},z \in \mathbb{Z}, t \in \mathbb{T}
\end{split}
\end{equation}
\begin{equation}
\begin{split}
\label{H2gen_shut_down}
N_{k,z} - n_{k,z,t} \geq \sum\limits_{\tau = t-\tau^\text{DOWN}_{k,z}}^{t} n^{\text{DOWN}}_{k,z,t}
\qquad \forall  k \in \mathbb{K},z \in \mathbb{Z}, t \in \mathbb{T}
\end{split}
\end{equation}

\subsection{Hydrogen Storage Constraints}
The cumulative H$_2$ in storage should be kept under its maximum capacity and above a volume for cushion gas requirement, as follows:
\begin{equation}
\begin{split}
V^{\text{STO}}_{s,z} \geq \sum\limits^{t}_{\tau = t_{0}}\left( h^{\text{CHA}}_{s,z,\tau}  \eta^{\text{STO}}_{s,z} - \dfrac{h^{\text{DIS}}_{s,z,\tau}}{\eta^{\text{STO}}_{s,z}} \right)\Delta t  \geq \underline{\text{R}}^{\text{STO}}_{s,z}  V^{\text{STO}}_{s,z}
\\\qquad \forall  z \in \mathbb{Z}, s \in \mathbb{S}, t \in \mathbb{T}
\end{split}
\end{equation}
where $ V^{\text{STO}}_{s,z} $ is the available storage capacity; $ \eta^{\text{STO}}_{s,z} $ is the charging and discharging efficiency (typically 100\% for H$_2$ storage); and $ \underline{\text{R}}^{\text{STO}}_{s,z} $ denotes the minimum ratio of stored H$_2$ to the maximum capacity. The charging rate of H$_2$ storage is physically capped by compression capability, as follows:
\begin{equation}
\begin{split}
H^{\text{STO}}_{s,z} \geq h^{\text{CHA}}_{s,z,t} \geq 0
\qquad \forall  z \in \mathbb{Z}, s \in \mathbb{S}, t \in \mathbb{T}
\end{split}
\end{equation}
% \begin{equation}
% \begin{split}
% h^{\text{DIS}}_{s,z,t} \geq 0
% \qquad \forall  z \in \mathbb{Z}, s \in \mathbb{S}, t \in \mathbb{T}
% \end{split}
% \end{equation}
% \subsection{Hydrogen Demand Constraints}
% The total hydrogen demand $ \text{D}_{z,t} $ is equal to the sum of hydrogen demands from transportation sector $ H2T_{z,t} $, heat sector $ H2H_{z,t} $, and industrial sector $ H2I_{z,t} $.
% \begin{equation}
% \begin{split}
% \label{}
% \text{D}_{z,t} = H2T_{z,t} + H2H_{z,t} + H2I_{z,t}
% \qquad \forall z \in \mathbb{Z}, t \in \mathbb{T}
% \end{split}
% \end{equation}
% The hydrogen demands from heat and industrial sectors are assumed to be flexible.
% \begin{equation}
% \begin{split}
% \label{}
% H2H^{min}_{z,t} \leq H2H_{z,t}  \leq H2H^{max}_{z,t}
% \qquad \forall z \in \mathbb{Z}, t \in \mathbb{T}
% \end{split}
% \end{equation}
% \begin{equation}
% \begin{split}
% \label{}
% H2I^{min}_{z,t} \leq H2I_{z,t}  \leq H2I^{max}_{z,t}
% \qquad \forall z \in \mathbb{Z}, t \in \mathbb{T}
% \end{split}
% \end{equation}
\subsection{Hydrogen Transmission Constraints}
% \subsubsection{Flexibility Assumptions}
% \begin{itemize}
% 	\item Pipeline, Truck can serve as storage via line packing and travel-time delays respectively. (Modeled)
% 	\item The filling time of truck/trailer is generally on the order of a few hours (Modeled indirectly through aboveground storage at point of charging/discharging)
% 	\item A truck can serve multiple routes (Currently NOT Modeled)
% 	\item boil-off losses from liquid hydrogen and liquefaction facility CAPEX in the case of liquid H2 based truck transport (Modeled)
% \end{itemize}
\subsubsection{Balance}
The total amount of H$_2$ transport between two zones is equal to the sum of the amounts through pipeline and truck. $ h^{\text{PIP}}_{z \rightarrow z',i,t} $ denotes the H$_2$ exchange at zone $ z $ through the pipeline $ i $ between zone $ z $ and $ z' $ during time $ t $. Positive $ h^{\text{PIP}}_{z \rightarrow z',i,t} $ represents that H$_2$ is delivered to zone $ z $, while negative value represents that H$_2$ flows out from zone $ z $. $ h^{\text{TRU}}_{z,j,t} $ denotes the H$_2$ exchange of the truck type $ j $ at zone $ z $ during time $ t $. Positive $ h^{\text{TRU}}_{z,j,t} $ represents that the truck is discharging H$_2$ to zone $ z $, while negative value represents that the truck is charging H$_2$ from zone $ z $.
\begin{equation}
\begin{split}
h^{\text{TRA}}_{z,t} = \sum\limits_{z' \in Z}\sum\limits_{i \in \mathbb{I}} h^{\text{PIP}}_{z \rightarrow z',i,t}
+
\sum\limits_{j \in \mathbb{J}} h^{\text{TRU}}_{z,j,t}
\qquad \forall  z \in \mathbb{Z}, t \in \mathbb{T}
\end{split}
\end{equation}

\subsubsection{Pipeline}
The H$_2$ exchange at zone $ z $ through the pipeline $ i $ between zone $ z $ and $ z' $ can be split into the H$_2$ delivering ($h^{\text{PIP+}}_{z \rightarrow z',i,t}$),  and the H$_2$ flowing out ($h^{\text{PIP-}}_{z \rightarrow z',i,t}$), as follows:
\begin{equation}
\begin{split}
h^{\text{PIP}}_{z \rightarrow z',i,t} =  h^{\text{PIP+}}_{z \rightarrow z',i,t} - h^{\text{PIP-}}_{z \rightarrow z',i,t} 
\quad \forall  z \rightarrow z^{'} \in \mathbb{B}, i \in \mathbb{I}, t \in \mathbb{T}
\end{split}
\end{equation}
The flow rate of H$_2$ through pipeline type $ i $ is capped by the operational limits of the pipeline $ i $, $ \overline{\text{F}}_{i} $, multiplied by the number of constructed pipeline $ i $ ($l_{z \rightarrow z',i}$), as follows:
\begin{equation}
\begin{split}
\overline{\text{F}}_{i}  l_{z \rightarrow z',i} \geq h^{\text{PIP+}}_{z \rightarrow z',i,t},h^{\text{PIP-}}_{z \rightarrow z',i,t} \geq 0
\quad \forall  z \rightarrow z^{'} \in \mathbb{B}, i \in \mathbb{I}, t \in \mathbb{T}
\end{split}
\end{equation}
% \begin{equation}
% \begin{split}
% \overline{\text{F}}_{i}  l_{z \rightarrow z',i} \geq h^{\text{PIP-}}_{z \rightarrow z',i,t} \geq 0
% \qquad \forall  z \rightarrow z^{'} \in \mathbb{B}, i \in \mathbb{I}, t \in \mathbb{T}
% \end{split}
% \end{equation}
The pipeline also has storage capability via line packing \cite{Clegg16}, modelled as follows:
\begin{equation}
\begin{split}
\overline{\text{E}}_{i}^{\text{PIP}}  l_{z \rightarrow z',i} \geq  &- \sum\limits^{t}_{\tau = t_{0}}\left( h^{\text{PIP}}_{z' \rightarrow z,i,\tau} + h^{\text{PIP}}_{z \rightarrow z',i,\tau} \right) \Delta t  \geq \underline{\text{R}}^{\text{PIP}}_{i}  \overline{\text{E}}_{i}^{\text{PIP}}  l_{z \rightarrow z',i}
\\\qquad &\forall  z' \in \mathbb{Z},z \in \mathbb{Z}, i \in \mathbb{I}, t \in \mathbb{T}
\end{split}
\end{equation}
% \begin{equation}
% \begin{split}
% - \sum\limits^{t}_{\tau = t_{0}}\left( h^{\text{PIP}}_{z' \rightarrow z,i,\tau} + h^{\text{PIP}}_{z \rightarrow z',i,\tau} \right) \Delta t  \leq \overline{\text{E}}_{i}^{\text{PIP}}  l_{z \rightarrow z',i}
% \\\qquad \forall  z' \in \mathbb{Z},z \in \mathbb{Z}, i \in \mathbb{I}, t \in \mathbb{T}
% \end{split}
% \end{equation}
where the maximum amount of H$_2$ that can be stored in pipeline is denoted by $ \overline{\text{E}}_{i}^{\text{PIP}} $, and the minimum amount, typically 0, is denoted by $ \underline{\text{R}}^{\text{PIP}}_{i} $.

\subsubsection{Truck}
\begin{figure*}[tb]
		\centering
		\includegraphics[width=0.7\textwidth]{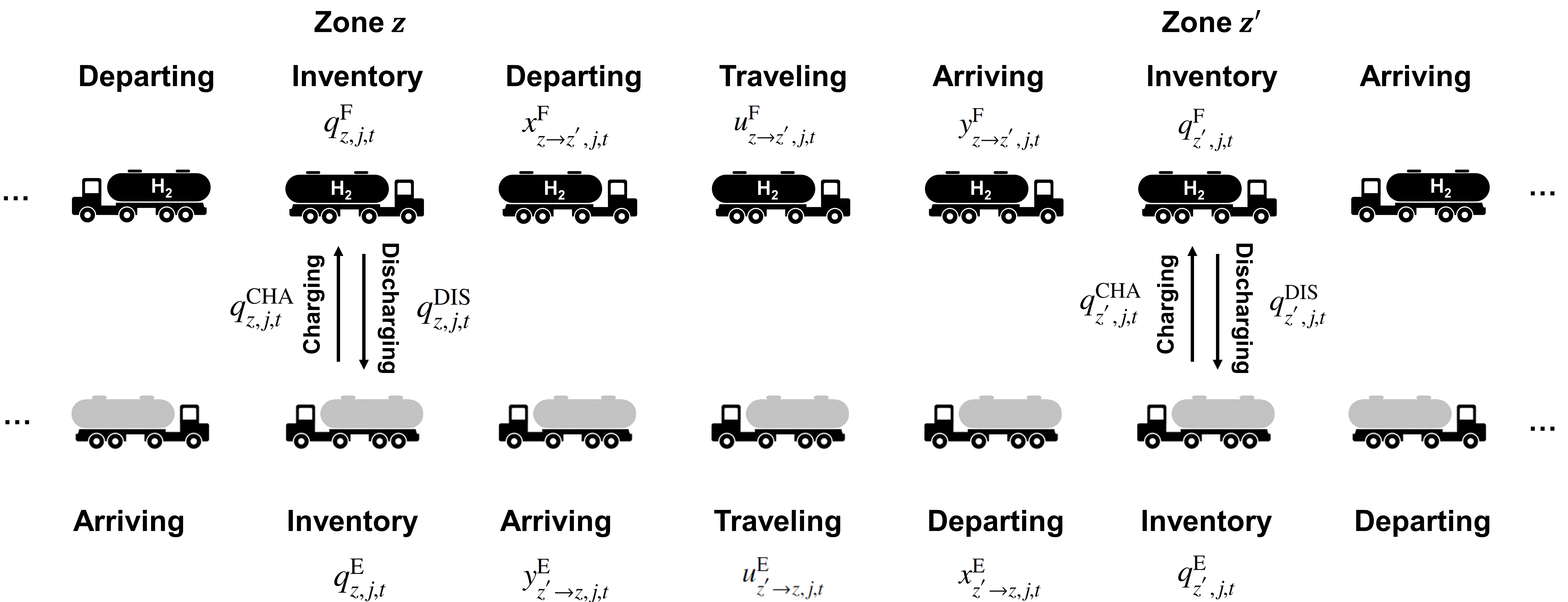}
		\caption{Schematic of the truck scheduling model. The upper row represents the trucks full of H$_2$, and the lower row represents empty trucks. Without loss of generality, zone $z$ is shown as the H$_2$ exporting zone, while zone $z'$ is the importing zone. Each truck can be dispatched to any zone and route.}
		\label{fig-truck}
\end{figure*}

We incorporate a flexible truck scheduling and routing model, which accurately captures the travelling delay of trucks and allows trucks to be shared across different routes and zones. We denote full and empty trucks with different sets of variables. Within each set, the state of trucks are further categorized into trucks in inventory at each zone and trucks in transit between each pair of zones, the latter category including departing, traveling, and arriving trucks. The schematic of the truck scheduling model is shown in Figure \ref{fig-truck}. While the number of trucks are integer in nature, we relax them to continuous variables to improve the computational tractability. In the case studies, we validate that this approximation is at little cost of modelling and decision accuracy.

The sum of full and empty trucks should equal the total number of invested trucks, as (\ref{full_empty_truck}). The full (empty) trucks include full (empty) trucks in transit and staying at each zones, as (\ref{full_empty_truck_2}):
\begin{equation}
\label{full_empty_truck}
\begin{split}
v^{\text{F}}_{j,t} +v^{\text{E}}_{j,t} = V_{j} 
\qquad \forall j \in \mathbb{J}, t \in \mathbb{T}
\end{split}
\end{equation}
\begin{equation}
\label{full_empty_truck_2}
\left\{
\begin{aligned}
v^{\text{F}}_{j,t} &= \sum\limits_{z \rightarrow z^{'} \in \mathbb{B}} u^{\text{F}}_{z \rightarrow z',j,t} + \sum\limits_{z \in \mathbb{Z}} q^{\text{F}}_{z,j,t} 
\\
v^{\text{E}}_{j,t} &= \sum\limits_{z \rightarrow z^{'} \in \mathbb{B}} u^{\text{E}}_{z \rightarrow z',j,t} + \sum\limits_{z \in \mathbb{Z}} q^{\text{E}}_{z,j,t}
\qquad \forall j \in \mathbb{J}, t \in \mathbb{T}
\end{aligned}
\right. 
\end{equation}

The change of the total number of full (empty) available trucks at zone $ z $ should equal the number of charged (discharged) trucks minus the number of discharged (charged) trucks at zone $ z $ plus the number of full (empty) trucks that just arrived minus the number of full (empty) trucks that just departed:
\begin{equation}
\left\{
\begin{aligned}
q^{\text{F}}_{z,j,t} - q^{\text{F}}_{z,j,t-1} 
= & q^{\text{CHA}}_{z ,j,t} - q^{\text{DIS}}_{z ,j,t} 
\\+& \sum\limits_{z' \in \mathbb{Z}} \left( - x^{\text{F}}_{z \rightarrow z',j,t-1} + y^{\text{F}}_{z \rightarrow z',j,t-1} \right)
\\
q^{\text{E}}_{z,j,t} - q^{\text{E}}_{z,j,t-1} 
=& - q^{\text{CHA}}_{z ,j,t} + q^{\text{DIS}}_{z ,j,t} 
\\+& \sum\limits_{z' \in \mathbb{Z}} \left( - x^{\text{E}}_{z \rightarrow z',j,t-1} + y^{\text{E}}_{z \rightarrow z',j,t-1} \right)
\\
\qquad \forall& z \in \mathbb{Z}, j \in \mathbb{J}, t \in \mathbb{T}
\end{aligned}
\right. 
\end{equation}

% \begin{equation}
% \begin{split}
% q^{\text{F}}_{z,j,t} - q^{\text{F}}_{z,j,t-1} 
% =& q^{\text{CHA}}_{z ,j,t} - q^{\text{DIS}}_{z ,j,t} 
% \\+& \sum\limits_{z' \in \mathbb{Z}} \left( - x^{\text{F}}_{z \rightarrow z',j,t-1} + y^{\text{F}}_{z \rightarrow z',j,t-1} \right)
% \\&\qquad \forall z \in \mathbb{Z}, j \in \mathbb{J}, t \in \mathbb{T}
% \end{split}
% \end{equation}
% \begin{equation}
% \begin{split}
% q^{\text{E}}_{z,j,t} - q^{\text{E}}_{z,j,t-1} 
% =& - q^{\text{CHA}}_{z ,j,t} + q^{\text{DIS}}_{z ,j,t} 
% \\+& \sum\limits_{z' \in \mathbb{Z}} \left( - x^{\text{E}}_{z \rightarrow z',j,t-1} + y^{\text{E}}_{z \rightarrow z',j,t-1} \right)
% \\&\qquad \forall z \in \mathbb{Z}, j \in \mathbb{J}, t \in \mathbb{T}
% \end{split}
% \end{equation}

The change of the total number of full (empty) trucks in transit from zone $ z $ to zone $ z' $ should equal the number of full (empty) trucks that just departed from zone $ z $ minus the number of full (empty) trucks that just arrived at zone $ z' $:
\begin{equation}
\left\{
\begin{aligned}
u^{\text{F}}_{z \rightarrow z',j,t} - u^{\text{F}}_{z \rightarrow z',j,t-1} 
= & x^{\text{F}}_{z \rightarrow z',j,t-1} - y^{\text{F}}_{z \rightarrow z',j,t-1}
\\
u^{\text{E}}_{z \rightarrow z',j,t} - u^{\text{E}}_{z \rightarrow z',j,t-1} 
= & x^{\text{E}}_{z \rightarrow z',j,t-1} - y^{\text{E}}_{z \rightarrow z',j,t-1}
\\
\qquad \forall& z \rightarrow z^{'} \in \mathbb{B}, j \in \mathbb{J}, t \in \mathbb{T}
\end{aligned}
\right. 
\end{equation}

% \begin{equation}
% \begin{split}
% u^{\text{F}}_{z \rightarrow z',j,t} - u^{\text{F}}_{z \rightarrow z',j,t-1} 
% = x^{\text{F}}_{z \rightarrow z',j,t-1} - y^{\text{F}}_{z \rightarrow z',j,t-1}
% \\\qquad \forall z \rightarrow z^{'} \in \mathbb{B}, j \in \mathbb{J}, t \in \mathbb{T}
% \end{split}
% \end{equation}
% \begin{equation}
% \begin{split}
% u^{\text{E}}_{z \rightarrow z',j,t} - u^{\text{E}}_{z \rightarrow z',j,t-1} 
% = x^{\text{E}}_{z \rightarrow z',j,t-1} - y^{\text{E}}_{z \rightarrow z',j,t-1}
% \\\qquad \forall z \rightarrow z^{'} \in \mathbb{B}, j \in \mathbb{J}, t \in \mathbb{T}
% \end{split}
% \end{equation}

The minimum travelling time delay is modelled in (\ref{Travel_delay_1})-(\ref{Travel_delay_2}). The number of full (empty) trucks in transit from $ z $ to $ z' $ at time $ t $ should be greater than the number of full (empty) truck departing $ z $ between time $ t-\Delta_{z \rightarrow z'}+1 $ and $ t $, as follows:
\begin{equation}
\label{Travel_delay_1}
\left\{
\begin{aligned}
u^{\text{F}}_{z \rightarrow z',j,t} \geq& \sum_{e=t-\Delta_{z \rightarrow z'}+1}^{e=t} x^{\text{F}}_{z \rightarrow z',j,e}  
\\
u^{\text{E}}_{z \rightarrow z',j,t} \geq& \sum_{e=t-\Delta_{z \rightarrow z'}+1}^{e=t} x^{\text{E}}_{z \rightarrow z',j,e}  
\qquad \forall z \rightarrow z' \in \mathbb{B}, j \in \mathbb{J}, t \in \mathbb{T}
\end{aligned}
\right. 
\end{equation}

% \begin{equation}
% \label{Travel_delay_1}
% \begin{split}
% u^{\text{F}}_{z \rightarrow z',j,t} \geq \sum_{e=t-\Delta_{z \rightarrow z'}+1}^{e=t} x^{\text{F}}_{z \rightarrow z',j,e}  
% \\ \qquad \forall z \rightarrow z' \in \mathbb{B}, j \in \mathbb{J}, t \in \mathbb{T}
% \end{split}
% \end{equation}
% \begin{equation}
% \begin{split}
% u^{\text{E}}_{z \rightarrow z',j,t} \geq \sum_{e=t-\Delta_{z \rightarrow z'}+1}^{e=t} x^{\text{E}}_{z \rightarrow z',j,e}  
% \\ \qquad \forall z \rightarrow z' \in \mathbb{B}, j \in \mathbb{J}, t \in \mathbb{T}
% \end{split}
% \end{equation}

The number of full (empty) trucks in transit from $ z $ to $ z' $ at time t should be greater than the number of full (empty) trucks arriving at $ z' $ from $ z $ between time $ t+1 $ and $ t+\Delta_{z \rightarrow z'} $, as follows:
\begin{equation}
\label{Travel_delay_2}
\left\{
\begin{aligned}
u^{\text{F}}_{z \rightarrow z',j,t}  \geq& \sum_{e=t+1}^{e=t+\Delta_{z \rightarrow z'}} y^{\text{F}}_{z \rightarrow z',j,e}
\\
u^{\text{E}}_{z \rightarrow z',j,t}  \geq& \sum_{e=t+1}^{e=t+\Delta_{z \rightarrow z'}} y^{\text{E}}_{z \rightarrow z',j,e} 
\qquad \forall z \rightarrow z^{'} \in \mathbb{B}, j \in \mathbb{J}, t \in \mathbb{T}
\end{aligned}
\right. 
\end{equation}

% \begin{equation}
% u^{\text{F}}_{z \rightarrow z',j,t}  \geq \sum_{e=t+1}^{e=t+\Delta_{z \rightarrow z'}} y^{\text{F}}_{z \rightarrow z',j,e}  \qquad \forall z \rightarrow z^{'} \in \mathbb{B}, j \in \mathbb{J}, t \in \mathbb{T}
% \end{equation}
% \begin{equation}
% \label{Travel_delay_2}
% u^{\text{E}}_{z \rightarrow z',j,t}  \geq \sum_{e=t+1}^{e=t+\Delta_{z \rightarrow z'}} y^{\text{E}}_{z \rightarrow z',j,e}  \qquad \forall z \rightarrow z^{'} \in \mathbb{B}, j \in \mathbb{J}, t \in \mathbb{T}
% \end{equation}

The amount of H$ _2 $ delivered to zone $ z $ should equal the truck capacity times the number of discharged trucks minus the number of charged trucks, adjusted by the H$_2$ boil-off loss during truck transportation and compression, as follows:
\begin{equation}
\begin{split}
h^{\text{TRU}}_{z,j,t}   = \left[(1 - \sigma_{j})q^{\text{DIS}}_{z ,j,t} - q^{\text{CHA}}_{z ,j,t}\right]  \overline{\text{E}}_{j}^{\text{TRU}}
\\\qquad \forall z \rightarrow z^{'} \in \mathbb{B}, j \in \mathbb{J}, t \in \mathbb{T}
\end{split}
\end{equation}

The charging capability of truck stations is limited by their compression or liquefaction capacity, as follows:
\begin{equation}
\begin{split}
q^{\text{CHA}}_{z,j,t}  \overline{\text{E}}_{j}^{\text{TRU}} \leq H^{\text{TRU}}_{z,j}
\qquad \forall z \in \mathbb{Z}, j \in \mathbb{J}, t \in \mathbb{T}
\end{split}
\end{equation}

%% file: Sections/4_results.tex
\section{Case Studies in the U.S. Northeast}
\label{case_study}
We apply the HSC planning model to evaluate the least-cost infrastructure needed to serve potential H$_2$ demand in transportation in the U.S. Northeast, based on simulating annual operations through operations over 20 representative weeks. The model is solved with Gurobi on a 40-core Intel Xeon Gold 6248 with 8 GB RAM.

\subfile{Sections/4.1_data_description}

\subsection{Results}
\subsubsection{Optimal H$_2$ generation mix}
\begin{figure}[htb]
		\centering
		\includegraphics[width=0.45\textwidth]{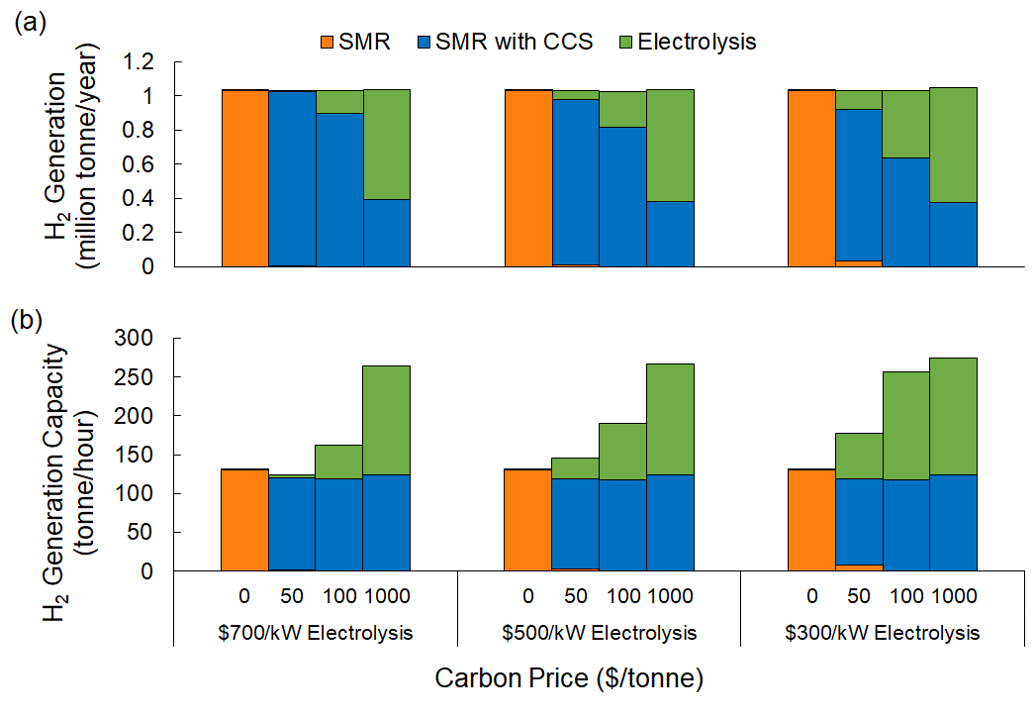}
		\caption{Optimal H$_2$ generation mix in the HSC under various electrolyzer cost and carbon price scenarios. (a) Annual H$_2$ generation (b) Invested H$_2$ generation capacity.}
		\label{fig-1}
\end{figure}

\begin{figure}[tb]
		\centering
		\includegraphics[width=0.45\textwidth]{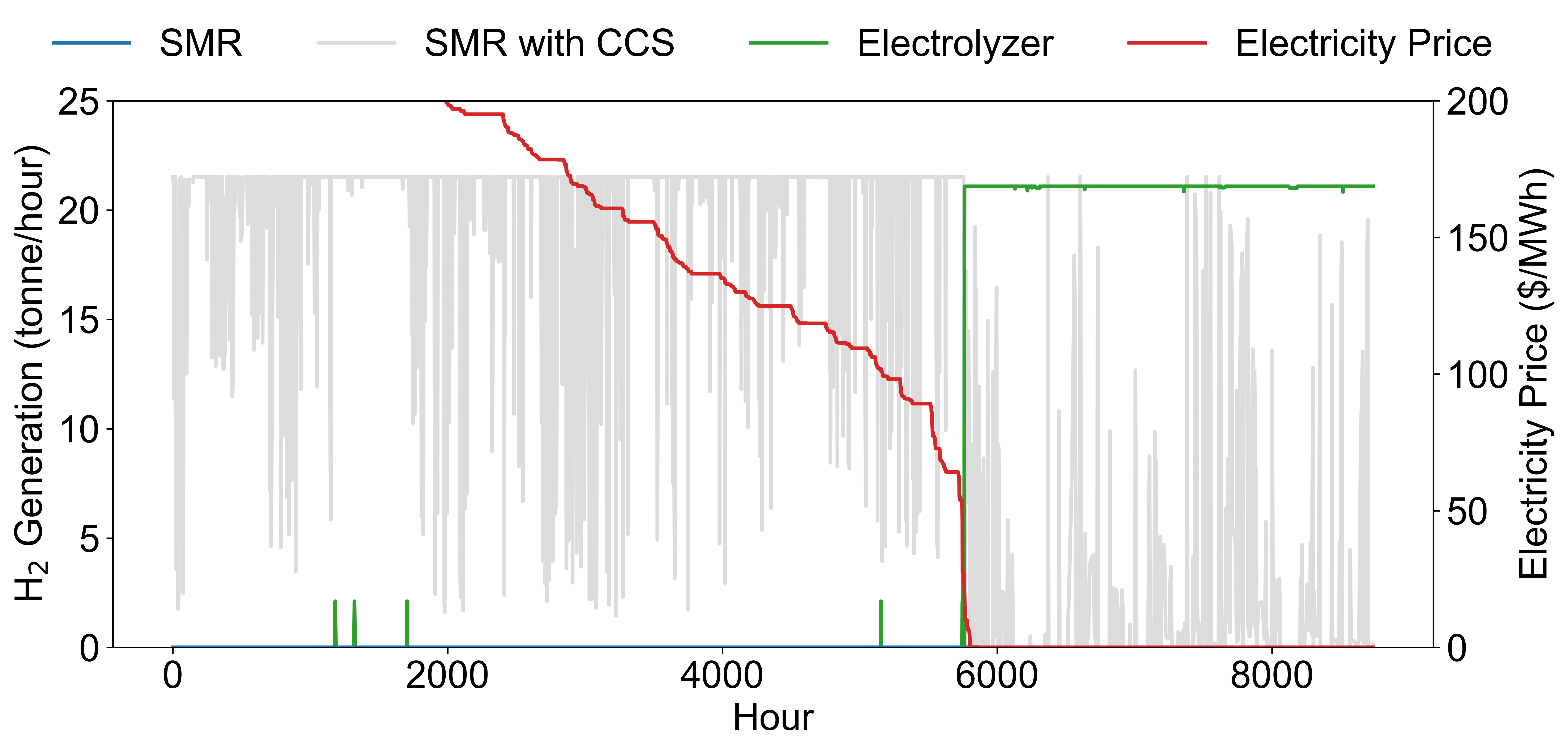}
		\caption{H$_2$ generation at different electricity prices at Zone 6 in the base case assuming \$100/tonne carbon price and \$300/kW electrolyzer cost.}
		\label{fig-2}
\end{figure}

Fig. \ref{fig-1} shows the optimal generation mixes in HSC with different carbon prices and unit capital cost of electrolyzer, which highlight several key observations. First, Fig. \ref{fig-1} (a) shows that as the carbon price increases, the H$_2$ generation switches from central SMR, to SMR with CCS, and then to electrolyzer. Second, with decreasing capital cost, electrolyzer is cost-effective for deployment at lower carbon prices. For example, as seen in Fig. \ref{fig-1} (b), installed electrolyzer capacity at \$100/tonne CO$_2$ price is 44 tonne/hour (2.4 GW) for an electrolyzer capital cost of \$700/kW as compared to 59 tonne/hour (3.1 GW) for \$300/kW. Thirdly, although carbon price 
increase favors an increasing share of generation from electrolyzer, it has a relatively small impact on the installed capacity of SMR with CCS in H$_2$ generation. SMR with CCS remains a cost-effective source of H$_2$ supply for time periods when electricity prices are very high for all carbon prices (see Fig. 3 (b)). The last observation can also be confirmed in Fig. \ref{fig-2}, where the electrolyzer produces H$_2$ when the electricity price is below approximately \$30/MWh. The results point to cost-effectiveness of electrolytic H$_2$ supply with moderate carbon policy (\$50/tonne or greater) and electrolyzer capital cost reduction (\$500/kW or lower). In the following discussions, we choose \$300/kW electrolyser cost and \$100/tonne carbon price as the benchmark scenario to explore other aspects of the HSC model outcomes.

\subsubsection{Truck scheduling model comparison}

TABLE \ref{tab:model_compare} compares the computational complexity and optimal investment decisions of our proposed model (linear truck scheduling) with that of the integer model (integer truck scheduling) and the model in existing literature ("Existing" model in TABLE \ref{tab:model_compare}, assuming a fixed transport capacity limit and ignoring (\ref{full_empty_truck})-(\ref{Travel_delay_2}) in the proposed model) \cite{almansoori16,kim17,li20,ochoa18,reu17,welder_18}. Two cases are examined, one is the base case, the other is a mini case with one representative week and gas truck only. The mini case highlights how relaxing truck scheduling variables to continuous rather than integer significantly reduces computational time while achieving almost the same planning outcome (only less than 0.3\% deviation in the H$_2$ truck and generation capacities). Notably, the integer model for the U.S. Northeast with 20 weeks is computationally intractable for the imposed model time limit of 24 h. The existing model without flexible truck scheduling over-invests both trucks (by 86\%) and stationary storage (by 213\%) because it does not allow trucks to be shared across different routes (as transmission) and zones (as storage) and to serve as storage simultaneously, resulting in an 11\% increase in H$_2$ supply cost. Moreover, electrolyzer generation increases (by 10\%) with truck sharing in the proposed model, which effectively enhances electrolyzer flexibility to respond to temporal variability in electricity prices while meeting H$_2$ demands. 

\begin{table}[tb]
  \begin{center}
  \caption{Sizes (after Presolve), Computation Times, and Outcomes among the Proposed Model, the Integer Model, and the Existing Model}
    \label{tab:model_compare}%
  \resizebox{0.48\textwidth}{!}
  {
    \begin{tabular}{lccccc}
    \toprule
      & \multicolumn{3}{c}{Base case} & \multicolumn{2}{c}{Mini case} \\
\cmidrule{2-6}      & \multicolumn{1}{l}{Proposed} & \multicolumn{1}{l}{Existing} & \multicolumn{1}{l}{Integer} & \multicolumn{1}{l}{Proposed} & \multicolumn{1}{l}{Integer} \\
\midrule
Continuous variable & 2064013 & 896205 &       & 67057 & 36842 \\
Integer variable & 0     & 0     &       & 0     & 47048 \\
Computation time & 1132 s & 138 s & > 24h & 7 s   & 2765 s \\
Unit hydrogen cost (\$/kg) & 2.011 & 2.23  & -     & 2.524 & 2.525 \\
Truck capacity (tonne) & 1544  & 2870  & -     & 1159  & 1157 \\
SMR (tonne/hour) & 117   & 118   & -     & 123   & 122 \\
Electrolyser (tonne/hour) & 139   & 126   & -     & 0     & 0 \\
Stationary storage (tonne) & 486   & 1520  & -     & 176   & 185 \\
\bottomrule
\end{tabular}%

    }
    \end{center}
\end{table}%

\subsubsection{Storage and transport functions of resources in the HSC}
In the base case, H$_2$ pipeline is not economically competitive. However, it could become attractive if the existing natural gas pipeline infrastructure can be retrofitted at relatively low cost for transporting H$_2$. Therefore, we also evaluate the case with 50\% pipeline cost in addition to the base case. The supplied storage capacities and the annual amounts of stored H$_2$ and transported H$_2$ in both cases are summarized in TABLE \ref{tab:contribution}. In the base case, mobile storage resources including gas and liquid trucks supply more than 70\% of the total H$_2$ storage capacity and throughput. Gas truck has higher contributions than liquid truck in both storage and transport, partially due to the lower cost and energy consumption of gas compression than liquefaction. In the case with 50\% pipeline cost, pipelines are built to handle over half of the total transported H$_2$, reducing the need in truck for transport compared to the base case. As the truck capacities decrease, the need in stationary storage rises in turn, while the mobile storage resources still contribute approximately 60\% of the total stored H$_2$. Although pipelines present significant storage capacities through line packing, they are hardly used as storage, because they incur less operational costs than trucks in transport and should be used for transport as much as possible.

\begin{table}[tb]
  \begin{center}
  \caption{The Storage and Transport Contributions of Various Technologies in the Base Case and the Case with 50\% Pipeline Cost}%
  \resizebox{0.48\textwidth}{!}
  {
    \begin{tabular}{llccccc}
    \toprule
          &       & Stationary  & \multicolumn{2}{c}{Mobile Storage} & \multirow{2}[2]{*}{Pipeline} & \multirow{2}[2]{*}{Total} \\
          &       & Storage & Gas truck & Liquid truck &       &  \\
    \midrule
    \multicolumn{1}{l}{\multirow{6}[2]{*}{Base case}} & Storage Capacity & \multirow{2}[1]{*}{486} & \multirow{2}[1]{*}{768} & \multirow{2}[1]{*}{776} & \multirow{2}[1]{*}{0} & \multirow{2}[1]{*}{2029} \\
          & (tonne) &       &       &       &       &  \\
          & Stored H$_2$ & \multirow{2}[0]{*}{0.064} & \multirow{2}[0]{*}{0.134} & \multirow{2}[0]{*}{0.042} & \multirow{2}[0]{*}{0} & \multirow{2}[0]{*}{0.246} \\
          & (M tonne/year) &       &       &       &       &  \\
          & Transported H$_2$  & \multirow{2}[1]{*}{0} & \multirow{2}[1]{*}{0.307} & \multirow{2}[1]{*}{0.061} & \multirow{2}[1]{*}{0} & \multirow{2}[1]{*}{0.379} \\
          & (M tonne/year) &       &       &       &       &  \\
    \midrule
          & Storage Capacity & \multirow{2}[1]{*}{734} & \multirow{2}[1]{*}{543} & \multirow{2}[1]{*}{389} & \multirow{2}[1]{*}{488} & \multirow{2}[1]{*}{2153} \\
          & (tonne) &       &       &       &       &  \\
    \multicolumn{1}{l}{50\%} & Stored H$_2$ & \multirow{2}[0]{*}{0.093} & \multirow{2}[0]{*}{0.105} & \multirow{2}[0]{*}{0.019} & \multirow{2}[0]{*}{0.001} & \multirow{2}[0]{*}{0.217} \\
    \multicolumn{1}{l}{Pipeline Cost} & (M tonne/year) &       &       &       &       &  \\
          & Transported H$_2$  & \multirow{2}[1]{*}{0} & \multirow{2}[1]{*}{0.152} & \multirow{2}[1]{*}{0.026} & \multirow{2}[1]{*}{0.218} & \multirow{2}[1]{*}{0.401} \\
          & (M tonne/year) &       &       &       &       &  \\
    \bottomrule
    \end{tabular}%
    }
  \label{tab:contribution}
    \end{center}
% \end{table}
% \begin{table}[htb]
  \begin{center}
  \caption{HSC Capacity Comparisons between Systems with and without Significant Electrolytic H$_2$ Generation and Different Pipeline Costs}
    \label{tab:pipe_gen}%
  \resizebox{0.48\textwidth}{!}
  {
    \begin{tabular}{lcccc}
\toprule
      & \multicolumn{2}{c}{\$300/kWh Electrolyzer} & \multicolumn{2}{c}{\$500/kWh Electrolyzer} \\
      & \multicolumn{2}{c}{\$100/tonne CO$_2$ Price} & \multicolumn{2}{c}{\$50/tonne CO$_2$ Price} \\
Pipeline Cost Factor & 100\% & 50\%  & 100\% & 50\% \\
\midrule
Pipeline Flow Capacity (tonne/hour) & 0     & 30    & 0     & 63 \\
Truck Capacity (tonne) & 1544  & 931   & 1339  & 349 \\
Storage Capacity (tonne) & 486   & 734   & 204   & 769 \\
Electrolyzer Capacity (tonne/hour) & 139   & 126   & 27    & 5 \\
\bottomrule
\end{tabular}%
    }
  
    \end{center}
\end{table}

\begin{figure}[tb]
		\centering
		\includegraphics[width=0.45\textwidth]{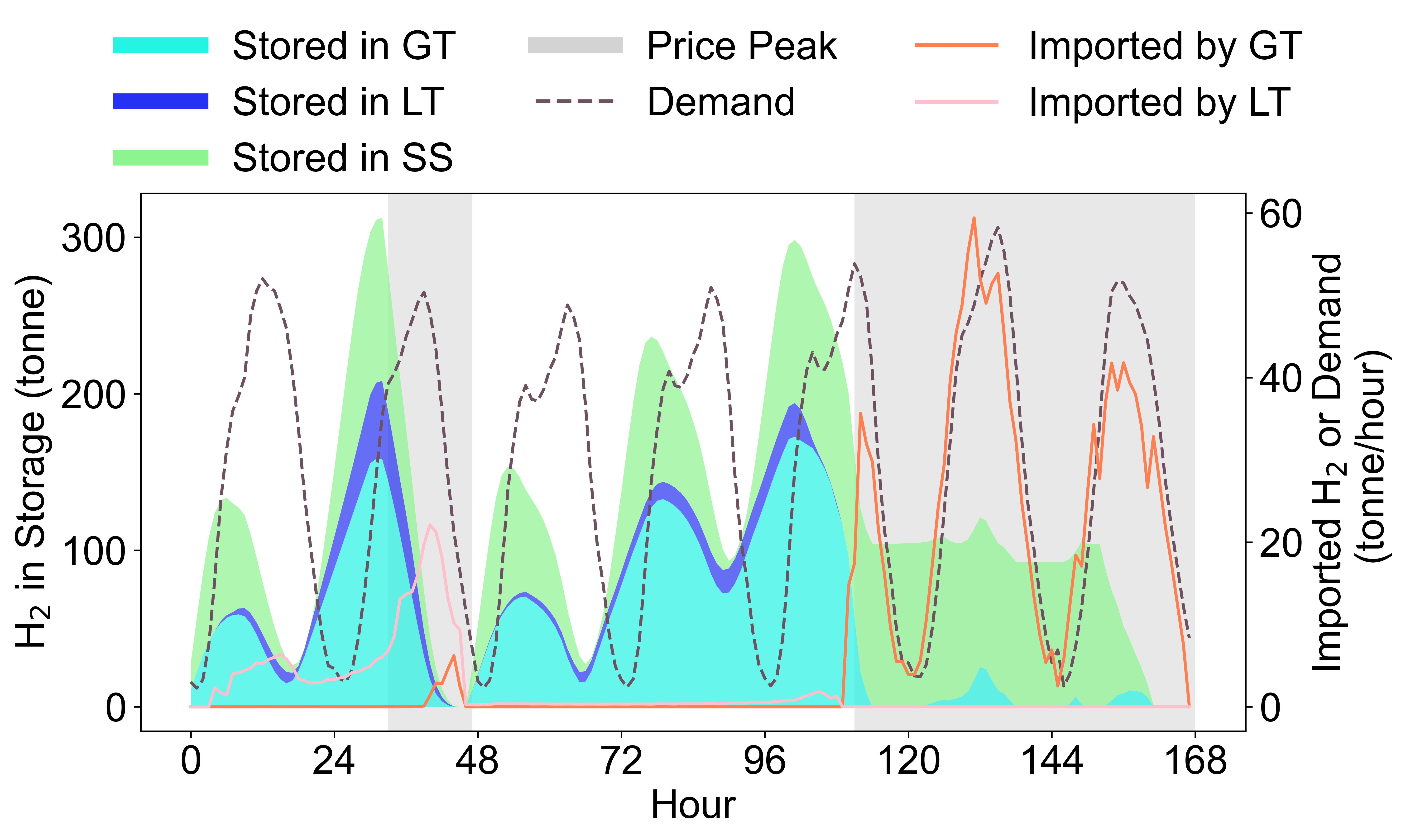}
		\caption{Profiles of stored and imported H$_2$ at Zone 4 in a representative week in the base case. GT: Gas Truck. LT: Liquid Truck. SS: Stationary Storage. 
%       (a) Base case (b) 50\% pipeline cost. 
        }
		\label{fig-3}
\end{figure}

Fig. \ref{fig-3} illustrates how mobile storage, both gas and liquid trucks, and stationary storage serve the HSC in response to the electricity price and demand fluctuations in a representative week. During off-peak electricity price periods, storage resources store H$_2$ during times of low H$_2$ demand and release H$_2$ during times of high H$_2$ demand. In the first short peak-price period (Hour 33-47), the trucks mostly serve as storage by releasing H$_2$. During the longer peak-price period (Hour 111-168), the storage resources are insufficient (or not cost-efficient to be over-built) to meet local demands, and the trucks are used for transporting and delivering H$_2$ from other zones installed with central SMR facilities.

TABLE \ref{tab:pipe_gen} compares the optimal infrastructures in systems with and without significant electrolytic H$_2$ generation. If pipeline is cheaper with retrofitting, more pipelines will be built in the system with less electrolytic H$_2$ generation (\$500/kWh electrolyzer, \$50/tonne CO$_2$ price) than the base case. This implies that trucks have synergy with distributed electrolysis, because they are more flexible and cost-efficient than pipelines to support intermittent electrolytic H$_2$ generation. Pipelines are more cost-efficient to cope with large and steady transmission demand in the centralized SMR production pathway, while electrolyzers can be deployed in a more distributed manner and thus complement 
smaller-scale and more flexible hydrogen transmission mechanisms like trucks.

% \begin{table}[tb]
%   \centering
%   \caption{Infrastructure Capacity Comparisons between Systems with and without Significant Electrolytic H$_2$ Generation and Different Pipeline Costs}
%     \label{tab:pipe_gen}%
%   \resizebox{0.48\textwidth}{!}
%   {
%     \begin{tabular}{lcccc}
% \toprule
%       & \multicolumn{2}{c}{\$300/kWh Electrolyzer} & \multicolumn{2}{c}{\$500/kWh Electrolyzer} \\
%       & \multicolumn{2}{c}{\$100/tonne CO$_2$ Price} & \multicolumn{2}{c}{\$50/tonne CO$_2$ Price} \\
% Pipeline Cost Factor & 100\% & 50\%  & 100\% & 50\% \\
% \midrule
% Pipeline Flow Capacity (tonne/hour) & 0     & 30    & 0     & 63 \\
% Truck Capacity (tonne) & 1544  & 931   & 1339  & 349 \\
% Storage Capacity (tonne) & 486   & 734   & 204   & 769 \\
% Electrolyzer Capacity (tonne/hour) & 139   & 126   & 27    & 5 \\
% \bottomrule
% \end{tabular}%
%     }

% \end{table}%

%% file: Sections/4.1_data_description.tex
\subsection{Data Descriptions}
Below is a brief description of the input parameters used to define the U.S. Northeast case studies, with further details available in the supporting document.

We focus on two types of H$_2$ generation technologies: electrolysis and natural gas fueled SMR (with and without CCS). We consider truck and pipeline as the key modes of H$_2$ transmission. We also model them as potential storage resources, in tandem with stationary H$_2$ storage. We model two types of trucks, based on handling H$_2$ as a cryogenic liquid or compressed gas, while the pipelines are considered as multiples of an 8" pipeline being built across different geographies. We do not consider geological H$_2$ storage as its availability in U.S Northeast region is uncertain \cite{ug_stor}.

\begin{figure}[tb]
    \centering
    \includegraphics[width = 0.48\textwidth]{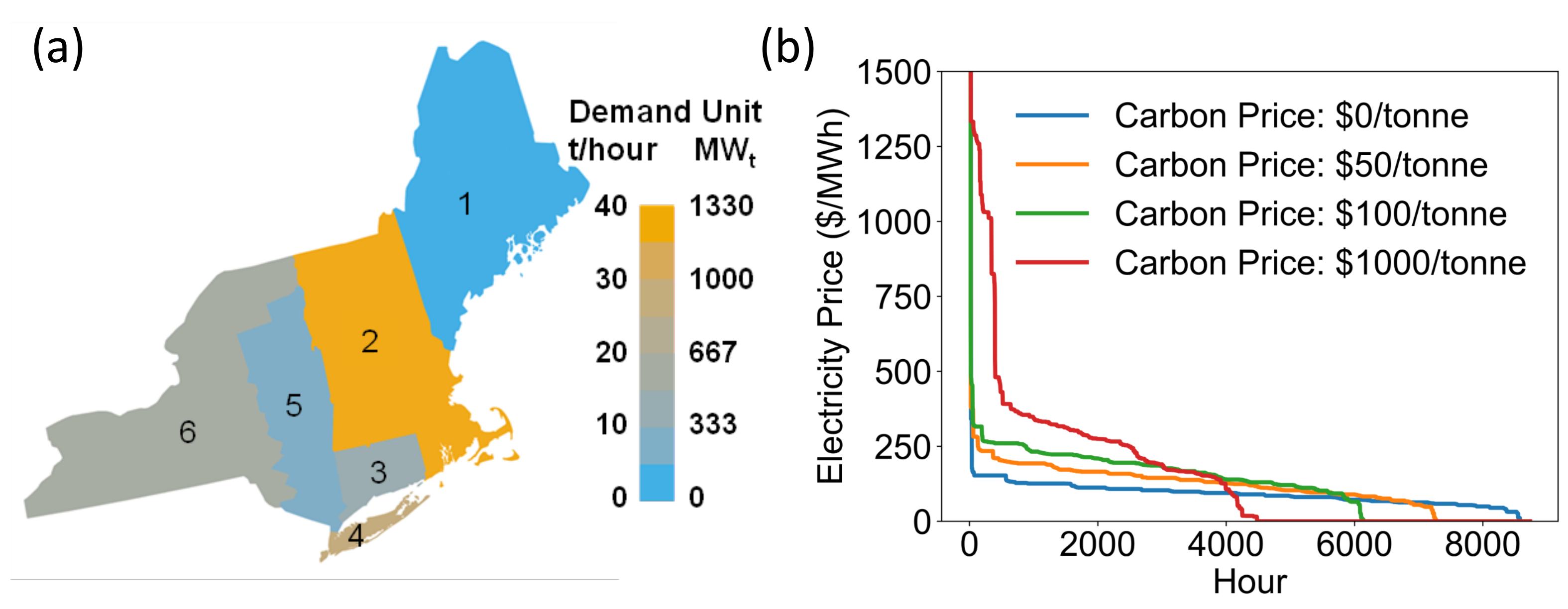}
    \caption{(a) Geographical zone classification for U.S. North-East and average H$_2$ demands for each zone. (b) Average price profiles generated from a power system capacity planning model for the studied U.S. Northeastern region.}
    \label{fig:demand_price}
\end{figure}

The U.S. Northeast region is represented in the model as 6 zones as shown in Fig. \ref{fig:demand_price} (a). The distances are measured by the distances between the polygon centroids of each zone. As Zone 2 and 4 are heavily urbanized, we do not allow central SMR to be built there in the model. The H$_2$ demands for each zone are developed based on available fuel consumption data and hourly charging profiles for both light- and heavy-duty fuel cell electric vehicles (FCEV) and the relative penetration of FCEV \cite{nexsnt,faf}. The average H$_{2}$ demands in the units of tonne/hour and MW$_\text{t}$ (lower heating value) are labeled in Fig. \ref{fig:demand_price} (a), based on a 20\% FCEV penetration.  

The electricity price profiles for the 20 representative weeks that are inputs to the HSC planning model are derived from a power systems capacity expansion model (CEM) for the U.S. Northeast under to the same carbon policy constraints as considered for the HSC model \cite{genx}.The CEM evaluates the least-cost mix of generation and storage technologies required to serve the prescribed electricity demand in 2050 based on modeling annual grid operations using the same representative weeks \cite{MALLAPRAGADA2020115390}. The electricity profiles under various carbon prices are presented in Fig. \ref{fig:demand_price} (b).

%\bibliography{Data_Sources}

%% file: Sections/5_conclusion.tex
\section{Conclusions}
\label{conclusion}
% DHARIK'S VERSION - 
Assessments on the role of hydrogen in future low-carbon energy systems have to contend with the multiple end-uses of hydrogen along with the plurality of technological options for its production, transport and storage. Here, we propose a least-cost HSC infrastructure planning model that accounts for a variety of operational and policy constraints, the model including: a) temporal and spatial variability in H$_2$ demand and renewable energy generation; b) a wide range of hydrogen-related technology options, and c) scalable representation of the flexibility offered by hydrogen truck and pipeline scheduling for transmission and storage functions. 

%
% YOUR VERSION - The multiple potential uses of hydrogen  in a future low-carbon energy system and the many technologies available for its production, transport and end-use requires considering interactions across the HSC but also between the HSC and other end uses. This paper proposes a HSC planning model that determines the least-cost mix of hydrogen production, storage, transmission, and compression facilities to meet hydrogen demands subject to a variety of operational and policy constraints, including variability in H$_2$ demand and availability and price of energy inputs (electricity). The model is highly configurable and capable of representing a wide range of conventional and emerging hydrogen-related technology options. We develop complex but scalable models for hydrogen truck and pipeline scheduling to fully capture their flexibility and their impacts on the HSC. The validation results indicate that modelling the trucks as mobile storage resource significantly changes the generation and transmission mix in HSC and thus is critical to the optimal planning decision making.

Application of the model to study HSC outcomes in the  U.S. Northeast reveals the potential for electrolyzer to be cost-effective under moderate carbon policy and electrolyzer cost reduction (\$50/tonne carbon price and \$500/kW unit capital cost). Enabling trucks to be shared across zones (as storage) and routes (as transmission) is critical to fully exploit their flexibility for VRE integration and HSC cost minimization.

If the hydrogen demand increases as it penetrates into the heating and industrial sectors in the future, the HSC will have a greater impact on the electricity price and becomes a price-maker. Then the planning of HSC requires coupled modelling of power systems and HSC, which is an area of future work. More complex and accurate modelling of pipeline flow rate and linepack and including the option of retrofitting natural gas pipeline are also interesting future directions.

%% file: Sections/Acknowledgment.tex
\section*{Acknowledgment}
We would like to thank Dr. Joe Powell, Dr. Mark Klokkenburg, and Dr. Robert Armstrong for their valuable advice on this work. We acknowledge the MIT SuperCloud and Lincoln Laboratory Supercomputing Center for providing resources that have contributed to the research results in this paper.